\documentclass[12pt,a4paper,leqno]{amsart}
\usepackage{amsmath}
\usepackage{amssymb}

\newtheorem{theorem}{Theorem}[section]
 
\newtheorem*{theoremA}{Theorem  A} 
\newtheorem*{theoremB}{Theorem  B} 
\newtheorem*{theoremC}{Theorem  C} 
 
\newtheorem{lemma}[theorem]{Lemma}
\newtheorem{proposition}[theorem]{Proposition}
\newtheorem{corollary}[theorem]{Corollary}

\newcommand{\R}{{\bf R}}       
\newcommand{\Z}{{\bf Z}}       

\catcode`\@=11 \@addtoreset{equation}{section} \catcode`\@=12

\title[ Nonlinear potentials and two weight trace inequalities] {Nonlinear
potentials and two weight trace inequalities for general dyadic  and 
 radial kernels} 
\author[Cascante]{Carme Cascante}
\address{Departament de Matem\`atica Aplicada i An\`alisi,
Facultat de Matem\`atiques, Universitat de Barcelona, Gran Via
585, 08071~Barcelona, Spain} \email{cascante@mat.ub.es}

\author[Ortega]{Joaquin M. Ortega }
\address{ Departament de Matem\`atica Aplicada i An\`alisi,
Facultat de Matem\`atiques,  Universitat de Barcelona, Gran Via
585, 08071~Barcelona, Spain} \email{ortega@mat.ub.es}

\author[Verbitsky]{Igor E. Verbitsky }
\address{Department of Mathematics, University of Missouri, 
Columbia, MO 65211, USA}
\email{igor@math.missouri.edu}
\thanks{First two authors partially supported by DGICYT Grant
BFM2002-04072-C02-01, and CIRIT Grant 2001SGR00172.}
\thanks{Third author partially supported
by NSF Grant DMS-0070623.}\subjclass{31C45, 46E35}

\keywords{Nonlinear potentials, Wolff's
inequality, two weight inequalities}

\begin{document}
\begin{abstract} 
We  study trace inequalities of the type
$$
\| T_k f\|_{L^q(d\mu)}\leq C \, \|f\|_{L^p(d\sigma)}, \qquad f \in
L^p(d\sigma),
$$  
in the ``upper triangle case''
$1 \leq q<p$  
for integral operators $T_k$ with positive kernels, where 
$d\sigma$ and $d\mu$  are   positive Borel measures  on $\R^n$. 
Our main tool is a generalization of Th. Wolff's inequality which gives 
 two-sided estimates of the energy ${\mathcal E}_{k,\,\sigma}
[\mu]=\int_{\R^n} (T_k [\mu])^{p'} \, d \sigma$ through  the $L^1(d\mu)$-norm
of an appropriate nonlinear potential $W_{k, \, \sigma}[\mu]$ associated with
the kernel $k$ and  measures $d\mu$, $d \sigma$. We initially work with a
dyadic integral operator with  kernel $K_{\mathcal D}(x, y) =
\sum_{Q\in{\mathcal D}}  K(Q) \chi_Q(x) \, \chi_Q(y)$, where $\mathcal
D=\{Q\}$ is the family of  all dyadic  cubes in $\R^n$, 
and $K: \, {\mathcal D}\rightarrow
\R^+$. The corresponding continuous versions of Wolff's
inequality and  trace inequalities are derived from their dyadic counterparts.

\end{abstract} \maketitle

\section{Introduction}  \label{introduction}

 In the present paper we  are concerned with  integral inequalities of the
type
 $$
\| T_k f\|_{L^q(d\mu)}\leq C \, \|f\|_{L^p(d\sigma)}, \qquad f \in
L^p(d\sigma),
$$  
in the ``upper triangle case''
$1\leq q<p$. Here 
$d\sigma$ and $d\mu$  are  locally finite positive Borel measures  on $\R^n$,
and
$$T_k f (x) = \int_{\R^n} K(x, y) \, f(y) \, d \sigma (y), \qquad x \in
\R^n,$$ is an integral operator with nonnegative kernel $K(x, y) \ge 0$.

An important part of this work is
a related study of generalized nonlinear potentials 
used  originally by Hedberg and Wolff 
\cite{hedbergwolff} in the special case of Riesz kernels $k_\alpha(x,y)=
|x-y|^{\alpha-n}$ ($0<\alpha<n$)   and when 
$d\sigma$ is 
Lebesgue measure on $\R^n$. Then 
$T_{k_\alpha} f(x) = k_\alpha \star f$ is the usual (linear) Riesz potential 
of $d \nu = f \, dx$. For a positive locally finite measure $\nu$ on $\R^n$, 
the corresponding nonlinear 
potential  is defined  by:
 $$W_{k_\alpha} [\nu](x) = \int_0^{+\infty} \left( \frac {\nu(B(x,r))}
{r^{n-\alpha p}}\right)^{p'-1} \frac {dr}{r},$$
where $B(x,r)$ is a ball of radius $r$ centered at $x \in \R^n$, and 
$p'=p/(p-1)$, $1<p< +\infty$. It follows that  
the energy   ${\mathcal E}_{k_\alpha}[\nu] = \left \Vert\, k_\alpha
\star \nu\right \Vert^{p'}_{L^{p'}(dx)}$ satisfies the classical  Wolff's  
inequality \cite{hedbergwolff}:
$$C_1 \,\int_{\R^n} W_{k_\alpha} [\nu] \, d \nu \le {\mathcal
E}_{k_\alpha}[\nu] \le C_2 \,  \int_{\R^n} W_{k_\alpha} [\nu]  \, d \nu,$$ 
where the constants $C_1, \, C_2$ 
do not depend on  $\nu$. Here only the upper estimate is nontrivial.  
(A thorough discussion of Wolff's inequality 
and its applications is given in  \cite{adamshedberg}.)

The
corresponding trace inequality  
$$||T_{k_\alpha} f||_{L^q(d \mu)} \le C \,
||f||_{L^p(dx)}, \qquad f \in L^p(dx),$$ 
in the case  $q=1$ by duality is equivalent to $W_{k_\alpha} [\mu] \in L^1(d
\mu)$.  As was shown by the authors in
\cite{cascanteortegaverbitsky}, 
for  $1<q<p$,  the preceding inequality is characterized by the condition 
$W_{k_\alpha} [\mu] \in L^{\frac {q(p-1)}{p-q}} (d \mu)$.  (Some
extensions of this result for more general kernels and weights  can be found
in  \cite{verbitsky}, \cite{cascanteortega} and
\cite{cascanteortegaverbitsky2}.)  Another characterization of this 
inequality  for $q<p$ stated in  capacitary terms   was obtained earlier by
Maz'ya and Netrusov  \cite{mazya}, \cite{mazyanetrusov}. 

 We observe that in the well-studied case $p\leq q<+\infty$, the above 
inequalities  have been fully characterized in several different ways,   for 
Riesz potentials and  Lebesgue  measure $dx$
in place of  $d\sigma$, in   \cite{adams}, \cite{kermansawyer}, 
 \cite{mazya}, \cite{mazyaverbitsky}, \cite{verbitsky}, and 
for a wide class
of integral operators $T_k$ in \cite{nazarovtreilvolberg}, 
\cite{sawyerwheeden}, \cite{verbitskywheeden}.

Our main goal is to unify and extend earlier results on Wolff's
inequality and trace inequalities for  more general two weight estimates
and  integral operators $T_k$ with  dyadic and radial nonincreasing kernels in
the case $q<p$.

We 
 first concentrate on  dyadic  integral operators
$T_{K_{\mathcal D}}$ introduced below, and the corresponding integral
inequalities. Later on we will show how continuous versions follow from their
dyadic analogues.    
 Let $ {\mathcal D}=\{ Q\}$ be the family of all dyadic cubes $Q$ in $\R^n$,
and $K: \, {\mathcal D}\rightarrow \R^+$. The kernel $K_{\mathcal D}(x,y)$ on
$\R^n\times\R^n$ is defined by  $$K_{\mathcal D}(x,y) = \sum_{Q\in{\mathcal
D}}  K(Q) \,  \chi_Q(x) \, \chi_Q(y),$$ 
where $\chi_Q$ is the characteristic function of $Q \in {\mathcal
D}$.

Let $\nu$ be a locally finite
positive Borel measure on $\R^n$, and let $f \in L^1_{\text{loc}} (d \nu)$. 
We define the  dyadic integral 
operator:
  $$ T_{K_{\mathcal D}}[f \, d\nu] (x) =
 \int_{\R^n} K_{\mathcal D} (x, y) f(y) \,
d\nu(y)=\sum_{Q\in{\mathcal D}}  K(Q)
\chi_Q(x) \int_Q f \, d \nu.$$   
 In case $f\equiv 1$, we set 
$$T_{K_{\mathcal D}}[\nu] (x) = \sum_{Q\in{\mathcal D}}  K(Q) \, \nu(Q) \, 
\chi_Q(x).$$

 If $0<q, \, p<+\infty$, and $\sigma$ and $\mu$ are locally finite Borel
measures on $\R^n$,  the corresponding dyadic trace inequality is given by:
\begin{equation}\label{equation0.1}
\int_{\R^n} \left \vert \, T_{K_{\mathcal D}}[f \, d\sigma] \, \right \vert^q
\,d\mu\leq  C \, \|f\|_{L^p(d\sigma)}^q,  \qquad f \in
L^p(d\sigma).
\end{equation}

Assume for a moment that $q, \, p> 1$. Duality then gives that 
 (\ref{equation0.1}) is equivalent to the inequality:
\begin{equation}\label{equation0.1dual}
\int_{\R^n} \left \vert \, T_{K_{\mathcal D}}[gd\mu]\right \vert^{p'} \,
d\sigma  \leq C \, \| g\|_{L^{q'}(d\mu)}^{p'}, \qquad g \in
L^{q'} (d\mu).
\end{equation}

The quantity on the left-hand side of (\ref{equation0.1dual}) is a generalized 
version of the discrete energy of $d \nu = gd\mu$.
For  positive locally finite Borel measures  $\nu$ and $\sigma$  on $\R^n$,
the discrete energy associated with $\nu$ and $\sigma$ is defined by   
(cf. \cite{hedbergwolff}):
\begin{equation}\label{discretenergy} 
{\mathcal E}_{K,\,\sigma}^{\mathcal
D}[\nu]= \int_{\R^n}  \left(T_{K_{\mathcal D}}[\nu]\right)^{p'}
\, d\sigma = \int_{\R^n} \left( \sum_{Q\in{\mathcal
D}}K(Q) \, \nu(Q)   \, \chi_Q(x)\right)^{p'}d\sigma(x). 
\end{equation} 
Fubini's theorem gives an alternative expression for ${\mathcal
E}_{K,\,\sigma}^{\mathcal D}$:
 $${\mathcal
E}_{K,\,\sigma}^{\mathcal D}[\nu]=\int_{\R^n}T_{K_{\mathcal D}}[(T_{K_{\mathcal
D}}[\nu])^{p'-1}d\sigma] \, d\nu,$$  
where  $T_{K_{\mathcal
D}}[(T_{K_{\mathcal D}}[\nu])^{p'-1}d\sigma]$ is a dyadic analogue of the 
nonlinear potential of Havin--Maz'ya originally defined for $d \sigma =dx$ 
and with $k_\alpha$ in place of $K_{\mathcal D}$ (see \cite{adamshedberg}, 
\cite{mazya}).

In the   special case where $d \sigma$  is Lebesgue measure
on $\R^n$,  $K(Q)=\frac 1 {|Q|^{1-(\alpha/n)}}$, $0<\alpha<n$ and 
$|Q|$ is the Lebesgue measure of $Q$, i.e., when  $K_{\mathcal D}$  is
a discrete Riesz kernel on $\R^n$, Hedberg and Wolff 
introduced a dyadic   nonlinear potential   defined by: 
$$ {\mathcal W}_{\alpha, \, dx}^{\mathcal D} [\nu](x)=
\sum_{Q\in{\mathcal D}} \left(  \frac{\nu(Q)}{r_Q^{n-\alpha
p}}\right)^{p'-1}\chi_Q(x). $$  
(Here $r_Q$ denotes the
side length  of $Q$.) A dyadic version of Wolff's inequality established in
\cite{hedbergwolff}  shows that,  for 
$1<p<+\infty$,  $$ C_1 \, {\mathcal E}_{\alpha,\,dx}^{\mathcal D}[\nu] \le 
\int_{\R^n}{\mathcal W}_{\alpha, \, dx}^{\mathcal D} [\nu](x) 
\, d\nu(x) \le C_2 \, {\mathcal E}_{\alpha,\,dx}^{\mathcal D}[\nu]. $$ 
Consequently,  the trace inequality (\ref{equation0.1}) holds for $q=1$,
$1<p< +\infty$, and $d \sigma=dx$ if and only if ${\mathcal W}_{\alpha, \,
dx}^{\mathcal D} [\mu]$ is in $L^1(d\mu)$.  For  $1<q<p<+\infty$, as was
shown in \cite{cascanteortegaverbitsky},  (\ref{equation0.1}) holds if and
only if ${\mathcal W}_{\alpha, \, dx}^{\mathcal D} [\mu] \in
L^{\frac{q(p-1)}{p-q}}(d\mu)$.

In this line of argument,  we need to define a suitable nonlinear  potential
associated with a pair of measures  $\nu$, $\sigma$ and  the kernel
$K_{\mathcal D}$ so that   it is  applicable to  characterization of the
trace inequality  (\ref{equation0.1}) for general dyadic kernels.   

Let $\nu$ and $\sigma$ be positive
locally finite Borel measures on $\R^n$. We denote by $\overline{K}(Q)$
the function  
$$ \overline{K}(Q)(x)=\frac1{\sigma(Q)}\sum_{Q'\subset Q} K(Q')
\sigma(Q')\chi_{Q'}(x). 
$$
For $x\in\R^n$, we set 
\begin{equation}\label{formula00.1}
{\mathcal{W}}_{K, \,\sigma}^{\mathcal D} [\nu](x)=\sum_{Q\in{\mathcal
D}} K(Q) \, \sigma(Q) \left( \int_{ Q} \overline{K}(Q)(y)d\nu(y) 
\right)^{p'-1}\chi_Q(x). \end{equation}

In \cite{cascanteortegaverbitsky2}, 
the following Wolff-type inequality was established for an arbitrary 
positive measure $\nu$ on $\R^n$, and $dx$ in place of $d \sigma$:
  $$
C_1 \, {\mathcal E}_{K,\,dx}^{\mathcal D} [\nu] \le 
\int_{\R^n}{\mathcal{W}}_{K, \, dx}^{\mathcal D} [\nu] \, d\nu \le C_2 \, 
{\mathcal E}_{K,\,dx}^{\mathcal D}[\nu], 
$$
where $C_1, \, C_2$ are constants which do not depend on $\nu$ and $K$.  We
note that some statements  in  \cite{cascanteortegaverbitsky2} (in 
particular, Lemma 2.6,
and consequently  Lemma 2.7,   Propositions 4.2 and 4.3) required certain
corrections outlined  in the erratum to that paper. They  are presented with
full proofs   in Section \ref{section01} below. Our main results  extend  the
preceding estimates, as well as  characterizations of the discrete trace
inequality  in \cite{cascanteortegaverbitsky},
\cite{cascanteortegaverbitsky2},  to arbitrary measures $\sigma, \, \mu$ on
$\R^n$.  In particular, we will prove the following theorems.

\begin{theoremA}\label{theoremA}
Let $K:{\mathcal D}\rightarrow\R^+$. Let $1<p<+\infty$, 
 and let $\nu$ and $\sigma$ be locally finite positive Borel
measures on $\R^n$. If  ${\mathcal E}_{K,\, \sigma}^{\mathcal
D}[\nu]$ and ${\mathcal{W}}_{K, \,  \sigma}^{\mathcal D} [\nu]$ are defined
respectively  by {\rm(\ref{discretenergy})} and  {\rm(\ref{formula00.1})},
then  $$
C_1 \, {\mathcal E}_{K,\, \sigma}^{\mathcal D}[\nu] \le 
\int_{\R^n}{\mathcal{W}}_{K, \,  \sigma}^{\mathcal D} [\nu] (x) 
\, d\nu (x)
\le C_2 \,  {\mathcal E}_{K,\, \sigma}^{\mathcal D}[\nu], 
$$
where $C_1, \, C_2$ are constants which do not depend on $\nu$ and $\sigma$.
\end{theoremA}

\begin{theoremB}\label{theoremB}
Let $K:{\mathcal D}\rightarrow\R^+$. Let $1\le q<p<+\infty$, 
 and let $\mu$ and $\sigma$ be locally finite positive Borel
measures on $\R^n$. 

 {\rm(i)} Suppose there exists a constant $C>0$ 
such that the trace inequality
$$\int_{\R^n} \left \vert \, T_{K_{\mathcal D}}[fd\sigma] \, \right \vert^q(x)
\,d\mu(x)\leq  C \, \|f\|_{L^p(d\sigma)}^q,  \qquad f\in
L^p(d\sigma),$$
 holds. Then 
 ${\mathcal W}_{K, \,\sigma}^{\mathcal D} [\mu]\in
L^{\frac{q(p-1)}{p-q}}(d\mu)$.

{\rm(ii)} Conversely, if 
${\mathcal W}_{K, \,\sigma}^{\mathcal D} [\mu]\in
L^{\frac{q(p-1)}{p-q}}(d\mu)$ then  the preceding trace inequality holds
provided the pair $(K, \sigma)$  satisfies  the dyadic logarithmic bounded
oscillation condition {\rm (\text{DLBO})}:
  $$\sup_{x \in Q} \overline{K}(Q)(x) \le A \, \inf_{x \in Q}
\overline{K}(Q)(x),$$ where $A$ does not depend on $Q \in {\mathcal D}$.  
 \end{theoremB} 
 
If $q=1$ then statement (ii) holds without the restriction $
(K, \sigma)\in\text{DLBO}$. (In this case Theorem B is, by duality, 
an immediate consequence of Theorem
A.)

In Section \ref{section02}, we obtain continuous analogues of Theorems A and B.
Here we state only a continuous version of the trace inequality for
convolution operators with radial  kernels $k(x) = k(|x|)$, 
$$T_k[f](x) = \int_{\R^n} k(|x-y|) \, f(y) \, d \sigma(y).$$
Here 
$k=k(r)$, $r>0$, is an arbitrary lower
semicontinuous  nonincreasing positive function.

The corresponding nonlinear potential is defined by  
$$
 {\mathcal W}_{k,\,\sigma}[\mu](x)=\int_0^{+\infty} k(r) \, 
\sigma(B(x,r)) \, 
\left(\int_{B(x,r)}\overline{k}(r)(y) \, d\mu(y)\right)^{p'-1}  \,
\frac{dr}{r},  $$
where 
$$
\overline{k}(r)(x)=\frac1{\sigma(B(x,r))}\int_0^r
k(s) \, \sigma(B(x,s)) \, \frac{ds}{s},$$ 
for $x\in\R^n$, $r>0$.

\begin{theoremC}\label{theoremC}
 Let $1\le q<p<+\infty$, 
 and let $\mu$ and $\sigma$ be locally finite positive Borel
measures on $\R^n$. Assume that  $\sigma$ satisfies a  doubling condition, and 
the pair $(k, \sigma)$  has  the  logarithmic bounded oscillation property
{\rm (\text{LBO})}:  $$\sup_{y \in B(x,r)} \overline{k}(r)(y) \le A \,
\inf_{y \in B(x,r)} \overline{k}(r)(y),$$ where $A$ does not depend on
$x\in\R^n, r>0$.  Then the following conditions are equivalent:

 {\rm(i)} There exists a constant $C>0$ 
such that the trace inequality
$$\int_{\R^n} \left \vert \, T_{k}[fd\sigma] \, \right \vert^q(x)
\,d\mu(x)\leq  C \, \|f\|_{L^p(d\sigma)}^q,  \qquad f\in
L^p(d\sigma),$$
 holds.

{\rm(ii)}  
${\mathcal W}_{k, \,\sigma} [\mu]\in
L^{\frac{q(p-1)}{p-q}}(d\mu)$.
  \end{theoremC} \bigskip
 
\noindent {\bf Remark 1.} We observe that  $(k, \sigma)\in \,${\rm\text{LBO}} 
if $d \sigma = dx$ is Lebesgue measure, for arbitrary nonincreasing 
radial kernels $k(r)$, 
or if $k(r)=r^{\alpha-n}$ is a Riesz kernel and $\sigma$ satisfies 
a reverse doubling condition of order $\gamma>n -\alpha$. (See details in 
Section \ref{section02}.)\bigskip

\noindent {\bf Remark 2.} An example given in Sec.~3 demonstrates that 
Theorem C is no longer true  for the nonlinear 
potential defined  by 
$$
 \overline{{\mathcal W}}_{k,\,\sigma}[\mu](x)=\int_0^{+\infty}  \, 
\overline{k}(r)(x) \, \sigma(B(x,r)) \, \left(\int_{B(x,r)}\overline{k}(r)(y)
\, d\mu(y)\right)^{p'-1}  \, \frac{dr}{r},  $$
in place of  $ {\mathcal W}_{k,\,\sigma}[\mu]$, even when $d \sigma$ is
Lebesgue measure and  $\overline{k}(r) = r^{-n} \int_0^r k(s) \, s^{n-1} \,
ds$ depends  only on $r$.

A similar example in Sec.~2 shows that  Theorems A and B fail  if
one replaces $ {\mathcal W}_{K, \,\sigma}^{\mathcal D}[\mu]$ by 
$$\overline{{\mathcal W}}_{K, \,\sigma}^{\mathcal D}
[\mu](x)=\sum_{Q\in{\mathcal D}} 
\overline{K}(Q)(x) \, \sigma(Q) \left( \int_{Q} \overline{K}(Q)(y) \, 
d\mu(y)  \right)^{p'-1}\chi_Q(x),$$
    even if 
$d \sigma$ is Lebesgue measure and $\overline{K}(Q)(x)$ 
is constant on $Q$.
\bigskip

 We conclude the introduction with a remark on our notation: we will adopt the
usual convention of using the same letter for various ``absolute'' constants 
(which may depend on $q, \, p$ and $n$) 
whose values may change in each occurrence, and we will write $A\preceq
B$ if there exists an absolute constant $M$ such that $A\leq MB$.
 We will say that two quantities $A$ and $B$ are equivalent if both
$A\preceq B$ and $B\preceq A$, and in that case we will write
$A\simeq B$.
\section{Discrete Wolff-type  and   trace inequalities}
\label{section01}

Let $K:{\mathcal D}\rightarrow \R^+$ where $\R^+=[0, +\infty)$. 
We consider the kernel $K_{\mathcal D}$ given by
$$
K_{\mathcal D}(x,y)=\sum_{Q\in{\mathcal D}}K(Q)\chi_Q(x)\chi_Q(y),
$$
for $x$, $y$ in $\R^n$. If $\mu$ is a positive locally finite Borel measure 
on $\R^n$, we define the operator $T_{K_{\mathcal D}}$ by
$$T_{K_{\mathcal D}}[\mu](x)=\sum_{Q\in{\mathcal D}}K(Q)\mu(Q)\chi_Q(x).$$

Suppose $\sigma$ is a fixed positive locally finite measure on $\R^n$. For 
$d \mu=f \, d\sigma$, where $f$ is a nonnegative Borel measurable function, 
we will simply  write 
$$T_{K_{\mathcal D}}[f](x)=\sum_{Q\in{\mathcal
D}}K(Q)\chi_Q(x)\int_Qf(y)d\sigma(y).$$  

If $1<p<+\infty$, the discrete energy of $\mu$ is given by:
$${\mathcal E}_{K_{\mathcal D}, \, \sigma}[\mu]=\int_{\R^n}
\left( T_{K_{\mathcal
D}}[\mu](x)\right)^{p'}d\sigma(x)=\int_{\R^n}\left(\sum_{Q\in{\mathcal
D}} K(Q)\mu(Q)\chi_Q(x)\right)^{p'}d\sigma(x).$$

We now define a suitable nonlinear potential  which
generalizes the classical Hedberg-Wolff potential. For $1<p<+\infty$ and $Q \in
{\mathcal D}$, we first define the function  
$$
\overline{K}(Q)(x)=\frac1{\sigma(Q)}\sum_{Q'\subset Q}  K(Q')
\sigma(Q')\chi_{Q'}(x). 
$$ 
Note that $\overline{K}(Q)(x)$ is supported on $Q$. If $d\sigma=dx$ is 
Lebesgue measure on $\R^n$ and $K(Q)$ depends only on the size of $Q$, i.e.
there exists a nonincreasing function $k:(0,+\infty)\rightarrow \R^+$ such
that for any $Q\in{\mathcal D}$ $K(Q)=k(r_Q)$, it is easy to check that
$$\overline{K}(Q)(x)\simeq \frac1{r_Q^n}\int_0^{r_Q} k(t) t^{n-1}dt.$$  (See
the proof of part (i) of Proposition \ref{examples} below.) Next, we set  
\begin{equation}\label{formula0.1} {\mathcal{W}}_{K, \,\sigma}^{\mathcal D}
[\mu](x)= \sum_{Q\in{\mathcal D}}K(Q)\sigma(Q)\left(   \int_{ Q}
\overline{K}(Q)(y)d\mu(y)  \right)^{p'-1}\chi_Q(x), \quad x\in\R^n.
\end{equation} It is worthwhile to observe that several other natural 
alternatives to  ${\mathcal{W}}_{K,\,\sigma}^{\mathcal D} [\mu]$ discussed in
\cite{cascanteortegaverbitsky2} fail to satisfy the desired analogue of
Wolff's inequality. (See an example at the end of this section.)

We will also deal with dyadic ``shifted'' versions of the above potential
defined by \begin{equation}\label{formula0.1.0}
{\mathcal{W}}_{K, \,\sigma}^{{\mathcal D}_z}[\mu](x)=
\sum_{Q\in{\mathcal D}_z}K(Q)\sigma(Q)\left(   \int_{ Q}
\overline{K}(Q)(y)d\mu(y)  \right)^{p'-1}\chi_Q(x), \quad x\in\R^n,
\end{equation}
where ${\mathcal D}_z$ denotes the shifted dyadic lattice 
${\mathcal D}_z={\mathcal D}+z=\{Q+z\}_{Q\in{\mathcal D}}$.

We also introduce a dyadic  maximal function associated 
with $K^{\mathcal D}$:
$M_K^{\mathcal D} [\mu]$  given by 
$$M_K^{\mathcal D} [\mu](x)=\sup_{x\in
Q}\frac1{\sigma(Q)}\left(\sum_{Q'\subset Q}K(Q')\sigma(Q')\mu(Q')\right).$$ 

We recall that if $d\sigma=dx$  is Lebesgue measure on $\R^n$ and $K(Q)$ 
depends only on the size of $Q$, this  maximal function can be
rewritten as the dyadic  maximal function considered in
\cite{kermansawyer}, namely $$M_K^{\mathcal D} [\mu](x)=\sup_{x\in
Q}\overline{K}(Q)(x)\mu(Q).$$

Indeed, if $Q\in{\mathcal D}$, and $x\in Q$, then for any $l\geq0$, 
$2^{-l}Q$ is the unique cube in ${\mathcal D}$ satisfying $x\in 2^{-l}Q$ and
$r_{2^{-l}Q}=2^{-l}r_Q$. We then have
\begin{equation}\label{fractionalmaximalfunction} 
\frac1{r_Q^n}\sum_{Q'\subset Q}k(r_{Q'})r_{Q'}^n\mu(Q')=\frac{\mu(Q)}{r_Q^n}
\sum_{l\geq 0} k(\frac{r_Q}{2^l}) \left(\frac{r_Q}{2^l}\right)^n\simeq \mu(Q)
\overline{K}(Q)(x). \end{equation} Our first result can be viewed as a
discrete version of the Wolff inequality \cite{hedbergwolff}  and  
Kerman-Sawyer inequality \cite{kermansawyer}  for   general measures $\sigma$
and dyadic kernels $K_{\mathcal D}(x,y)$.

\begin{theorem}\label{theoremdiscrete} Let $K:{\mathcal D}\rightarrow\R^+$
and   $1<p<+\infty$. Let   $\mu$ and $\sigma$ be locally finite
positive Borel measures on $\R^n$. Then the  following
quantities are equivalent:
\begin{equation*}\begin{split}
& (a) \,\, {\mathcal E}_{K_{\mathcal D}, \,\sigma} [\mu]=  \int_{\R^n}\left(
\sum_{Q\in{\mathcal D}} K( Q) \mu(Q)\chi_Q(x) \right)^{p'}d\sigma(x);\\&
(b) \int_{\R^n}{\mathcal W}_{K, \, \sigma}^{\mathcal D} [\mu] \, d\mu =
\sum_{Q\in{\mathcal D}}K(Q)\sigma(Q)\mu(Q)
\left(\int_Q \overline{K}( Q)(y)d\mu(y)\right)^{p'-1};
\\& (c) \int_{\R^n}M_{K }^{\mathcal D} [\mu]^{p'} \, d\sigma.
\end{split}\end{equation*}

\end{theorem}
{\bf Proof of Theorem \ref{theoremdiscrete}:}\par
   
      Let $1<s<+\infty$,  $\Lambda=(\lambda_Q)_{Q\in{\mathcal D}}$, 
$\lambda_Q\in\R^+$,  and let $\sigma$ be a positive locally finite Borel
measure. Our standing  assumption will be that $\lambda_Q = 0$ if $\sigma(Q)
=0$. We will also  follow the  convention that $0 \cdot \infty = 0$.  

 We define  \begin{equation*}\begin{split}  & A_1(\Lambda)=
\int_{\R^n}\left( \sum_{Q\in{\mathcal D}} 
\frac{\lambda_Q}{\sigma(Q)}\chi_Q(x) \right)^sd\sigma(x),\\   &
A_2(\Lambda)=\sum_{Q\in{\mathcal D}} \lambda_Q \left(  
\frac{1}{\sigma(Q)}\sum_{Q'\subset Q}  \lambda_{Q'}\right)^{s-1},\\ 
&A_3(\Lambda)=\int_{\R^n}\sup_{x\in Q}\left( 
\frac1{\sigma(Q)}\sum_{Q'\subset Q}\lambda_{Q'}\right)^s  d\sigma(x).
\end{split}\end{equation*}   The proof of the theorem will be a consequence of
the following proposition. 

\begin{proposition}\label{proposition4.1}  Let 
$\sigma$ be a positive locally finite Borel  measure on $\R^n$.  Let  $1<s<
\infty$.  Then there exist constants $C_i >0, \, i = 1,2, 3$, which depend
only  on $s$,  such that, for any $\Lambda=(\lambda_Q)_{Q\in{\mathcal  D}}$,
$\lambda_Q\in \R^+$,  $$  A_1 (\Lambda) \leq C_1 \, A_2(\Lambda) \leq C_2 \, 
A_3(\Lambda)\leq C_3 \, A_1(\Lambda).$$      \end{proposition} 

 {\bf Proof of
Proposition \ref{proposition4.1}:}\par     We begin by observing that if $s\geq
1$, then for any
$\Lambda=(\lambda_Q)_{Q\in{\mathcal D}}$, $\lambda_Q\in\R^+$, and $x\in\R^n$,
we have   
\begin{equation}\label{sumationparts}   
\left(\sum_{Q\in{\mathcal
D}}\lambda_Q \, \chi_Q(x)\right)^s \leq s \,  \sum_{Q\in{\mathcal
D}}\lambda_Q \, \chi_Q(x) \left( \sum_{Q'\subset Q} \lambda_{Q'} \,
\chi_{Q'}(x) \right)^{s-1}.   
\end{equation}

We first prove (\ref{sumationparts}) under the assumption that 
$$\sum_{Q\in{\mathcal D}}\lambda_Q \, \chi_Q(x) < +\infty.
$$
Note that, for a fixed $x \in \R^n$, the 
 dyadic cubes containing $x$ form a nested family of cubes. 
Hence using the elementary inequality $b^s - a^s \le s \, (b-a) \, b^{s-1}$, 
$0\le a \le b$, $1 \le s < \infty$,  we
obtain:  $$ 
\left(\sum_{Q' \subset Q} \lambda_{Q'} \, \chi_{Q'} (x)\right)^s - 
\left(\sum_{Q' \subsetneqq Q} \lambda_{Q'} \, \chi_{Q'} (x)\right)^s \le \, 
s \, \lambda_Q \, \chi_Q(x) \left( \sum_{Q'\subset Q} \lambda_{Q'} \,
\chi_{Q'}(x) \right)^{s-1}. 
$$ 
From this (\ref{sumationparts}) follows by a telescoping sum argument,  
taking the sums of both sides  over all dyadic cubes $Q$ that contain $x$.

If $\sum_{Q\in{\mathcal
D}}\lambda_Q \, \chi_Q(x) = +\infty,$ but 
$\sum_{Q \subset Q_0} \lambda_{Q} \, \chi_{Q} (x)< +\infty$ for some (and
hence every) dyadic cube $Q_0$ which contains $x$ then  (\ref{sumationparts})
follows  by the same argument as above  taking the sums over all $Q \subset
Q_0$  and then letting $|Q_0| \to + \infty$. Finally, in the case where 
$\sum_{Q \subset Q_0}\lambda_Q \, \chi_Q(x) = +\infty$ for some 
$Q_0$, both sides of  (\ref{sumationparts}) are obviously 
infinite. This completes the proof of (\ref{sumationparts}).

        We now  prove $A_1(\Lambda)\leq C_1 \,  A_2(\Lambda)$  
for $1 < s\leq 2$ (this obviously holds for  $s=1$ and $C_1=1$ as well). 
We  may assume without loss of generality that there are only a finite number
of $\lambda_Q\neq 0$. By  (\ref{sumationparts}),  
\begin{equation*}\begin{split}  & A_1(\Lambda) = \int_{\R^n}\left(
\sum_{Q\in{\mathcal D}}  \frac{\lambda_Q}{\sigma(Q)}\chi_Q(x) \right)^s
d\sigma(x)\\ \leq &  s \,  \sum_{Q\in{\mathcal D}}
\frac{\lambda_Q}{\sigma(Q)}\int_Q  \left( \sum_{Q'\subset Q}  
\frac{\lambda_{Q'}}{\sigma(Q')}\chi_{Q'}(x)\right)^{s-1}d\sigma(x).  
\end{split}\end{equation*}  H\"older's inequality with exponent $\frac1{s-1}>
1$ gives:   \begin{equation*}\begin{split}  & \frac{1}{\sigma(Q)}\int_Q \left(
\sum_{Q'\subset Q}  
\frac{\lambda_{Q'}}{\sigma(Q')}\chi_{Q'}(x)\right)^{s-1}d\sigma(x)\\  \leq
&\left(  \frac{1}{\sigma(Q)}\int_Q  \sum_{Q'\subset Q}  
\frac{\lambda_{Q'}}{\sigma(Q')}\chi_{Q'}(x)d\sigma(x)\right)^{s-1}=  \left( 
\frac{1}{\sigma(Q)}\sum_{Q'\subset Q}   \lambda_{Q'}\right)^{s-1}. 
\end{split}\end{equation*}  Consequently, $$A_1(\Lambda)\leq s \, 
\sum_{Q\in{\mathcal D}}  \lambda_Q \left( \frac{1}{\sigma(Q)}\sum_{Q'\subset 
Q}\lambda_{Q'}\right)^{s-1} = s \,  A_2(\lambda),$$ which proves   the desired
inequality with $C_1=s$  for $1<s\leq 2$.  A similar inequality in the case
$s>2$ is proved by induction. For integer $k \ge 2$  we assume   that the
inequality  $A_1(\Lambda)\leq C_1 (s)\,  A_2(\Lambda)$  holds for any  $k-1< s
\leq k$, and  have to  show that it also holds for $k < s\le k+1$. By
(\ref{sumationparts}),    \begin{equation*}\begin{split}  & A_1(\Lambda) =
\int_{\R^n}\left( \sum_{Q\in{\mathcal D}} 
\frac{\lambda_Q}{\sigma(Q)}\chi_Q(x) \right)^s \, d\sigma(x)\\ \leq &  s \, 
\sum_{Q\in{\mathcal D}} \frac{\lambda_Q}{\sigma(Q)}\int_Q  \left(
\sum_{Q'\subset Q}   \frac{\lambda_{Q'}}{\sigma(Q')}\chi_{Q'}(x)\right)^{s-1}
\, d\sigma(x).   \end{split}\end{equation*}  Applying the induction hypothesis
for  $k-1<s-1\leq k$, with the  measure $\chi_Q \, \sigma$, and  the set
$(\lambda_{Q'})_{Q'}$,  where $\lambda_{Q'}=0$ for cubes $Q' \not\subset Q$,
we obtain:  \begin{equation*}\begin{split}  &
\sum_Q\frac{\lambda_Q}{\sigma(Q)} \int_Q  \left( \sum_{Q'\subset Q}  
\frac{\lambda_{Q'}}{\sigma(Q')}\chi_{Q'}(x)\right)^{s-1} \, d\sigma(x)\\ \leq
&  C_1(s-1) \,  \sum_Q\frac{\lambda_Q}{\sigma(Q)}\sum_{Q'\subset Q} 
\lambda_{Q'} \left( \frac{1}{\sigma(Q')}\sum_{Q''\subset 
Q'}\lambda_{Q''}\right)^{s-2} \\  =  & C_1(s-1) \sum_{Q'} \, \lambda_{Q'}
\left( \frac{1}{\sigma(Q')}\sum_{Q''\subset   Q'}\lambda_{Q''}\right)^{s-2}
\sum_{Q' \subset Q}\, \frac{\lambda_Q}{\sigma(Q)}  \\ \leq & C_1(s-1)
\int_{\R^n}  \sum_{Q'} \frac{\lambda_{Q'}}{\sigma(Q')} \, \chi_{Q'}(x) \,
\left( \frac{1}{\sigma(Q')}\sum_{Q''\subset  Q'}\lambda_{Q''}\right)^{s-2} \, 
\left (\sum_{Q}  \frac{\lambda_Q}{\sigma(Q)} \chi_Q (x)\right)  d \sigma(x).
\end{split}\end{equation*}  By H\"older's inequality for sums with exponents 
$s-1$  and  $(s-1)/(s-2)$ (note that $s-1>k-1\ge1$), we have
\begin{equation*}\begin{split} & \sum_{Q'} \frac{\lambda_{Q'}}{\sigma(Q')} \, 
\chi_{Q'}(x) \, \left( \frac{1}{\sigma(Q')}\sum_{Q''\subset 
Q'}\lambda_{Q''}\right)^{s-2}  \\ \leq &  
\left(\sum_{Q'}\frac{\lambda_{Q'}}{\sigma(Q')} \, 
\chi_{Q'}(x)\right)^{1/(s-1)}  \left(\sum_{Q'}\frac{\lambda_{Q'}}{\sigma(Q')}
\,  \chi_{Q'}(x)  \left( \frac{1}{\sigma(Q')}\sum_{Q''\subset 
Q'}\lambda_{Q''}\right)^{s-1} \right)^{(s-2)/(s-1)}.
\end{split}\end{equation*}  Substituting this estimate into the right-hand
side of the preceding inequality, we obtain:  \begin{equation*}\begin{split} 
& A_1(\Lambda) \leq  s \, C_1(s-1) \int_{\R^n} 
\left(\sum_{Q'}\frac{\lambda_{Q'}}{\sigma(Q')} \, 
\chi_{Q'}(x)\right)^{1/(s-1) +1}  \\ \times &\,
\left(\sum_{Q'}\frac{\lambda_{Q'}}{\sigma(Q')} \,    \chi_{Q'}(x) \left(
\frac{1}{\sigma(Q')}\sum_{Q''\subset  Q'}\lambda_{Q''}\right)^{s-1}
\right)^{(s-2)/(s-1)} \, d \sigma(x). \end{split} \end{equation*}  Applying
now H\"older's inequality for integrals with exponents  $s-1$   and 
$(s-1)/(s-2)$, we have:  \begin{equation*}\begin{split}  & A_1(\Lambda) \leq
s \, C_1 (s-1) \left(\int_{\R^n} 
\left(\sum_{Q'}\frac{\lambda_{Q'}}{\sigma(Q')} \,  \chi_{Q'}(x)\right)^{s} d
\sigma(x) \right)^{\frac1{s-1}}  \\ \times &
\left(\int_{\R^n}\sum_{Q'}\frac{\lambda_{Q'}}{\sigma(Q')} \,  \chi_{Q'}(x) 
\left( \frac{1}{\sigma(Q')}\sum_{Q''\subset  Q'}\lambda_{Q''}\right)^{s-1} \,
d \sigma(x)\right)^{\frac{s-2}{s-1}}\\ & = s \, C_1 (s-1)
A_1(\Lambda)^{\frac1{s-1} } A_2(\Lambda)^{\frac{s-2}{s-1}}.  
\end{split}\end{equation*}   From this it
follows that $ A_1(\Lambda) \leq C_1 (s)  \,  A_2(\Lambda)$ where $C_1(s) =
\left ( s \, C_1(s-1) \right )^\frac {s-1} {s-2}$,  for $k< s\leq 
k+1$, and hence for every $1<s<\infty$.

 Next, let $1<s<+\infty$ . Then  
\begin{equation*}\begin{split} 
 &A_3(\Lambda)=\int_{\R^n}\sup_{x\in Q}\left( 
\frac1{\sigma(Q)}\sum_{Q'\subset  Q}\lambda_{Q'}\right)^{s}d\sigma(x)\\ & \leq
\int_{\R^n}  M_\sigma^{H\!L,\,{\mathcal D}}\left( \sum_{Q\in{\mathcal D }} 
\frac{\lambda_Q}{\sigma(Q)}\chi_Q(x)\right)^sd\sigma(x) \leq  CA_1(\Lambda), 
\end{split}\end{equation*}  where the dyadic
Hardy-Littlewood maximal  function $M_\sigma^{H\!L, \, {\mathcal D}}[\nu]$ is defined by
$$M_\sigma^{H\!L,\, {\mathcal D}}[\nu] (x) =\sup_{x\in Q} \,
\frac{\nu(Q)}{\mu(Q)}.$$ 
Here we have used the fact that, for $d \nu = |f| \, d \sigma$, the operator 
$M_\sigma^{H\!L,\, {\mathcal D}}: f \to M_\sigma^{H\!L,\, {\mathcal D}}[f]$ is 
bounded on $L^s(\sigma)$ for $1<s \le +\infty$. This proves  the estimate
$A_3(\Lambda)\leq C  \, A_1(\Lambda)$ with $C =  ||M_\sigma^{H\!L,\,{\mathcal
D}}||^s_{L^s(\sigma)\to L^s(\sigma)}$.

 We now prove the  inequality 
$ A_1(\Lambda) \le C \, A_3(\Lambda)$
for
$1<s<\infty$. Using the estimate 
$A_1(\Lambda)\leq C_1 \, A_2(\Lambda)$ established above, we have:  
\begin{equation*}\begin{split}  &A_1(\Lambda)= \int_{\R^n} \left( \sum_{Q \in
{\mathcal D}}  \frac{\lambda_Q}{\sigma(Q)} \chi_Q(x)\right)^sd\sigma(x)\leq  
 C_1 \sum_{Q \in {\mathcal D}}\lambda_Q \left(  \frac1{\sigma(Q)}
\sum_{Q'\subset Q} \lambda_{Q'}\right)^{s-1}\\  & = C_1\int_{\R^n}\sum_{Q \in
{\mathcal D}}  \frac{\lambda_Q}{\sigma(Q)} \chi_Q(x) \left( 
\frac{1}{\sigma(Q)}\sum_{Q'\subset 
Q}\lambda_{Q'}\right)^{s-1}d\sigma(x)\\ & \leq C_1\int_{\R^n} \left( 
\sum_{Q\in{\mathcal D}} \frac{\lambda_Q}{\sigma(Q)}  \chi_Q(x)\right) \left(
\sup_{x\in Q}\frac1{\sigma(Q)}  \sum_{Q'\subset Q}
\lambda_{Q'}\right)^{s-1}d\sigma(x)\\ & \leq C_1 
A_1(\Lambda)^\frac1{s}A_3(\Lambda)^\frac1{s'},  
\end{split}\end{equation*} 
where in the last estimate we have used H\"older's inequality  with exponents
$s$ and $s'$. This chain of inequalities yields $$A_1(\Lambda)  \leq C_1 \,
A_1(\Lambda)^\frac1{s}A_3(\Lambda)^\frac1{s'},$$  which implies that  
$A_1(\Lambda) \leq C_1^{s'} \,A_3(\Lambda)$. It remains to show that  
$A_2(\Lambda)  \leq     C A_3(\Lambda).$  Note that
    $$A_2(\Lambda)
=\sum_{Q\in {\mathcal D}}\lambda_Q \left( \frac{1}{\sigma(Q)} \sum_{Q'\subset
Q}\lambda_{Q'} \right)^{s-1} =\sum_{Q\in {\mathcal D}}\lambda_Q\mu_Q,$$ where
$\mu_Q=\left(\frac{1}{\sigma(Q)} \sum_{Q'\subset Q}\lambda_{Q'}
\right)^{s-1}$. Then 
\begin{equation*}\begin{split} &\sum_{Q\in {\mathcal
D}}\lambda_Q\mu_Q= \int_{\R^n} \left( \sum_{Q\in{\mathcal D}} \frac{\lambda_Q
\mu_Q}{\sigma(Q)} \chi_Q(x)\right)d\sigma(x)\\ & \leq \int_{\R^n} \left(
\sum_{Q\in{\mathcal D}} \frac{\lambda_Q }{\sigma(Q)}
\chi_Q(x)\right)\sup_{x\in Q}\mu_Qd\sigma(x).
\end{split}\end{equation*} 
H\"older's inequality gives that the above is
bounded by \begin{equation*}\begin{split} &\left(\int_{\R^n} \left(
\sum_{Q\in{\mathcal D}} \frac{\lambda_Q }{\sigma(Q)}
\chi_Q(x)\right)^sd\sigma(x)\right)^\frac1{s}\left( \int_{\R^n}\sup_{x\in
Q}\mu_Q ^{s'}d\sigma(x)\right)^\frac1{s'}\\
&= \left(\int_{\R^n} \left( \sum_{Q\in{\mathcal D}} \frac{\lambda_Q
}{\sigma(Q)} \chi_Q(x)\right)^sd\sigma(x)\right)^\frac1{s} \left(
\int_{\R^n}\sup_{x\in Q}\left( \frac1{\sigma(Q)}\sum_{Q'\subset
Q}\lambda_{Q'}\right)^sd\sigma(x)  \right)^\frac1{s'}.
\end{split}\end{equation*} Hence $$A_2(\Lambda)=\sum_{Q\in{\mathcal D}} \lambda_Q\mu_Q \leq
A_1(\Lambda)^\frac1{s}A_3( \Lambda)^\frac1{s'}.$$    Since
$A_1(\Lambda)\leq C A_3(\Lambda)$, we finally obtain $A_2(\Lambda) \leq  C\,
A_3(\Lambda)$. The proof of the proposition is  complete.
\qed 

 Theorem \ref{theoremdiscrete} follows immediately from 
Proposition \ref{proposition4.1} with $\lambda_Q=K(Q ) 
\mu(Q)\sigma(Q)$     and $s=p'$. 

We are now in a position to characterize  
the trace inequality in the 
case $q=1$. 

\begin{corollary}\label{corollary0.1}
Let $K:{\mathcal D}\rightarrow\R^+$. Let $
\sigma, \mu$ be  locally finite positive Borel measures on $\R^n$, and let 
$1<p<+\infty$. Then the following statements are equivalent:

{\rm (i)} There exists $C>0$ such that for any $f\in L^p(d\sigma)$, $f \ge 0$, 
 $$\int_{\R^n}
T_{K_{\mathcal D}}[f] \, d\mu \leq C \, ||f||_{L^p(d\sigma)}.$$

 {\rm (ii)}
${\mathcal W}_{K,\,\sigma}^{\mathcal D} [\mu](x)=\sum_{Q\in{\mathcal
D}}\sigma(Q) \, K(Q) \, \chi_Q(x) \left( \int_Q
\overline{K}(Q)(y) \, d\mu(y)\right)^{p'-1}\in L^1(d\mu)$.
\end{corollary}
{\bf Proof of Corollary \ref{corollary0.1}:}\par
The trace inequality (i) can be restated equivalently via Fubini's theorem
as  $$
\int_{\R^n} T_{K_{\mathcal D}}[\mu]\, f \, d\sigma \leq
C \, \|f\|_{L^p(d\sigma)}, $$
which by duality is equivalent to  ${\mathcal E}_{K,\,\sigma}^{\mathcal
D}[\mu]<+\infty$. Now Theorem \ref{theoremdiscrete} gives that the 
dyadic energy is finite if and only if   ${\mathcal
W}_{K, \, \sigma}^{\mathcal D} [\mu]$ belongs to $L^1(d\mu)$.\qed    

 In what follows we will
restrict ourselves to functions $K:{\mathcal D}\rightarrow\R^+$ and measures
$\sigma$ satisfying an extra assumption, namely that 
for any $Q\in {\mathcal D}$,  there exists a constant that
we will denote by $\overline{K}(Q)$ such that 
\begin{equation}\label{extraassumption}  \frac1{C} \, \overline{K}(Q) \, 
\chi_Q(x) \leq \overline{K}(Q)(x)\leq C \, \overline{K}(Q) \, \chi_Q(x),
  \end{equation}  
where $C$ does not depend on $Q$. 
In
other words, the oscillation of the function $\ln \overline{K}(Q)(x)$ on $Q$
is bounded by a constant which is independent of $Q\in{\mathcal D}$. In this
case we will  say that the pair $(K,\sigma)$ satisfies the  {\em dyadic
logarithmic bounded oscillation\/}  property, or DLBO, and simply write
$(K,\sigma)\in \text{DLBO}$.  Observe that strictly speaking, 
$\overline{K}(Q)$ is not unequivocally defined, but since  it behaves like 
$\overline{K}(Q)(x)$ up to a multiplicative constant, we find this
notation appealing.   

 We  first give  some examples.    
\begin{proposition}\label{examples} {\rm (i)} Suppose that $d \sigma=dx$ is 
Lebesgue measure on $\R^n$ and $K(Q)$  depends only on the size  of 
$Q$, i.e., there exists a  nonincreasing function $k:(0,+\infty)\rightarrow
\R^+$ such that for any $Q\in{\mathcal D}$, $K(Q)=k(r_Q)$ where $r_Q$ is the
side length of  $Q$. Then the pair $(K,dx)$ satisfies property 
{\rm (\ref{extraassumption})} where 
$$\overline{K}(Q) = 
\frac1{r_Q^n}\int_0^{r_Q} k(t) \, t^{n-1} \, dt.$$        
{\rm (ii)} Let   
 $\displaystyle{K(Q)=\frac1{r_Q^{n-\alpha}}}$, $0<\alpha<n$, be  a discrete
Riesz kernel. Suppose  that    $\sigma$ is a dyadic reverse doubling measure: 
$\sigma \in DRD_\gamma$  for some $\gamma>n-\alpha$, i.e., there exists $C>0$
such that for any $j\geq0$, $Q\in {\mathcal D}$, 
\begin{equation}\label{mesuresdoblants1}   \sigma(2^jQ)\geq
C \, 2^{j\gamma} \, \sigma(Q),  \end{equation}  
where $2^jQ$ is the unique dyadic cube
in ${\mathcal D}$ such that $Q\subset 2^j Q$ and $r_{2^j Q}=2^j r_Q$.  We then
have that  $(K,\sigma)$ satisfies property {\rm \text{DLBO}} with
$\overline{K}(Q)=K(Q)$. Conversely, if $K$ is a discrete Riesz kernel and the
pair $(K,\sigma)$ satisfies property  {\rm \text{DLBO}} with
$\overline{K}(Q)=K(Q)$, then $\sigma$ is a $DRD_\gamma$ measure for some
$\gamma>n-\alpha$.

  \end{proposition}

 {\bf Proof of Proposition \ref{examples}:}\par
 We begin with (i). We observe that if $ Q\in{\mathcal D}$ and $x\in Q$, 
for any $l\geq 0$ 
there exists a unique $Q_l\subset Q$ in ${\mathcal D}$ such that $x\in Q_l$
and $r_{Q_l}=2^{-l}r_Q$. Thus   \begin{equation*}\begin{split}
  \overline{K}(Q)(x)  & =\frac1{r_Q^n} \sum_{Q'\subset Q} k(r_{Q'}) \,
r_{Q'}^n \,  \chi_{Q'}(x) \\& =
 \frac1{r_Q^n} \sum_{l\geq0} k(2^{-l}r_Q) \, (2^{-l}r_Q)^n \\& \leq \, C \, 
\frac1{r_Q^n}\int_0^{r_Q} k(t) \,  t^{n-1} \, dt.
\end{split}\end{equation*}

For the converse estimate, the fact that $k$ is nonincreasing gives
  \begin{equation*}\begin{split}
  \overline{K}(Q)(x)\geq 
 \frac1{r_Q^n} \sum_{l\geq1} k(2^{-l}r_Q) \, (2^{-l}r_Q)^n\geq  C 
\frac1{r_Q^n}\int_0^{ r_Q } k(t) \, t^{n-1} \, dt .
\end{split}\end{equation*}

Now we prove (ii). If $\sigma$ satisfies (\ref{mesuresdoblants1}),
$Q\in{\mathcal D}$, $x\in Q$ and $r_Q=2^{-k}$, and for any $l\geq0$, $2^{-l}Q$
is the unique cube in ${\mathcal D}$ such that $x\in 2^{-l}Q$ and
$r_{2^{-l}Q}=2^{-l}r_Q$, then   
\begin{equation*}\begin{split}
\overline{K}(Q)(x) &=\frac1{\sigma(Q)} \sum_{Q'\subset Q}
\frac{\sigma(Q')}{r_{Q'}^{n-\alpha}}\chi_{Q'}(x)\leq
\frac1{\sigma(Q)}\sum_{l\geq 0}   \sigma(2^{-l}Q) 2^{(l+k)(n-\alpha)}\\&  \leq
C  \sum_{l\geq 0}2^{ (l+k)(n-\alpha)}2^{-l\gamma}=   C
2^{k(n-\alpha)}\sum_{l\geq 0}   2^{ l(n-\alpha-\gamma)}\leq C
\frac{1}{r_Q^{n-\alpha}},  \end{split}\end{equation*}  since
$n-\alpha-\gamma<0$.  Obviously,  
$$\overline{K}(Q)(x)\geq 
\frac1{\sigma(Q)}\frac{\sigma(Q)}{r_Q^{n-\alpha}}
\chi_Q(x)=\frac1{r_Q^{n-\alpha}}.$$

 Suppose  now that (\ref{extraassumption}) with $\overline{K}(Q)=K(Q)$ holds.  
Let $Q\in{\mathcal D}$,   $x\in Q$ and $r_Q=2^{-k}$, where, as before, for
any $l\geq0$, $2^{-l}Q$ is the unique cube in ${\mathcal D}$ satisfying 
$x\in 2^{-l}Q$ and $r_{2^{-l}Q}=2^{-l}r_Q$. We then have: 
\begin{equation*}\begin{split}  C \, \frac{\sigma(Q)}{r_Q^{n-\alpha}} & \geq
\sum_{l\geq0} \frac{\sigma(2^{-l}Q)}{r_{2^{-l}Q}^{n-\alpha}}=  
\frac{\sigma(Q)}{r_Q^{n-\alpha}} + \sum_{l\geq1}
\frac{\sigma(2^{-l}Q)}{r_{2^{-l}Q}^{n-\alpha}} \\& \geq
\left(\frac1{C}+1\right)   \sum_{l\geq1}
\frac{\sigma(2^{-l}Q)}{r_{2^{-l}Q}^{n-\alpha}}\geq\cdots\geq \left(
\frac1{C}+1\right)^m \sum_{l\geq m}
\frac{\sigma(2^{-l}Q)}{r_{2^{-l}Q}^{n-\alpha}}\\& \geq \left(
\frac1{C}+1\right)^m   \frac{\sigma(2^{-m}Q)}{r_{2^{-m}Q}^{n-\alpha}}. 
\end{split}\end{equation*}    
Since $r_{2^{-m}Q}=2^{-m}r_Q$,  it follows that 
$$\sigma(Q)\geq
\frac1{C}  \left( (\frac1{C}+1) \, 2^{n-\alpha}\right)^m\sigma(2^{-m}Q),
$$  and
(\ref{mesuresdoblants1}) holds with $\gamma=n-\alpha+\log_2
\left(\frac{1}{C}+1 \right)$.

 \begin{theorem}\label{theorem4.2}
Let $K:{\mathcal D}\rightarrow\R^+$. Let $\mu$ and $\sigma$ 
be  locally finite positive Borel
measures on $\R^n$, $1<q<p<+\infty$. We then have:

{\rm(a)} If there exists $C>0$ such that for any $f\geq 0$,
\begin{equation}\label{trace}\left(\int_{\R^n}\left(
T_{K_{\mathcal D}}[f](x)\right)^q d\mu(x)\right)^\frac1{q}\leq C \, 
||f||_p,\end{equation} 
then $\displaystyle{{\mathcal W}_{K, \,\sigma}^{\mathcal D}
[\mu]=\sum_{Q\in{\mathcal D}}\sigma(Q) \, K(Q) \, \chi_Q\left(\int_Q
\overline{K}(Q)(y) \, d\mu(y)\right)^{p'-1} \in L^\frac{q(p-1)}{p-q}(d\mu)}$.

{\rm(b)} Conversely, suppose that in addition $(K,\sigma)\in${\rm\text{DLBO}}.
If  $\displaystyle{{\mathcal W}_{K, \,\sigma}^{\mathcal D}[\mu]\in
L^\frac{q(p-1)}{p-q}(d\mu)}$, then {\rm(\ref{trace})} holds.

\end{theorem}
{\bf Proof of Theorem \ref{theorem4.2}:}\par 
 
Duality gives  an equivalent reformulation of (\ref{trace}), namely,   
\begin{equation}||T_{K_{\mathcal D}} [gd\mu]||_{L^{p'}(d\sigma)}\leq C \, 
||g||_{L^{q'}(d\mu)}, \label{3.2}\end{equation} for any $g\in
L^{q'}(d\mu)$, $g\geq0$.

Theorem \ref{theoremdiscrete} applied to the positive measures
$gd\mu$ and $\sigma$, gives:
\begin{equation*}\begin{split} ||T_{K_{\mathcal D}} 
[gd\mu]||_{L^{p'}(d\sigma)}^{p'} & \geq  
C \int_{\R^n} {\mathcal W}_{K, \, \sigma}^{\mathcal D}
[gd\mu] \, g \,  d\mu \\ & =  C\sum_{Q\in{\mathcal D}} K(Q) \,\sigma(Q) \int_Q
g(x) \, d\mu(x)\left(\int_Q\overline{K}(Q)(x) \, g(x) \, d\mu(x)\right)^{p'-1}.
\end{split}\end{equation*}   
Assume that (\ref{trace}) (or  equivalently  
(\ref{3.2}) ) holds. We then have  that for any $g\in L^{q'}(d\mu)$, $g \ge
0$,  $$\sum_{Q\in{\mathcal D}}
K(Q)
\,
\sigma(Q)
\int_Q g(x) \, d\mu(x)\left(\int_Q\overline{K}(Q)(x)
\, g(x) \, d\mu(x)\right)^{p'-1} \leq C \, ||g||_{L^{q'}(d\mu)}^{p'}. $$  Let
$c_Q=K(Q) \, \sigma(Q) \, \mu(Q)\left(
\int_Q \overline{K}(Q)(x) \, d\mu(x)\right)^{p'-1}$. For $\psi\in
L^{\frac{q'}{p'}}(d\mu)$, $\psi\geq0$, let $$g(x)=\left( M_\mu^{H\!L,
\,{\mathcal D}}[\psi]\right)^\frac1{p'}(x):=\left( \sup_{x\in
Q}\frac1{\mu(Q)}\int_Q \psi(y) \, d\mu(y)\right)^\frac1{p'}.$$ The above
estimate  together with the $L^{q'}(d\mu)$-boundedness of the dyadic 
Hardy-Littlewood maximal function $M_\mu^{H\!L, \,{\mathcal D}}$, gives  
$$
\sum_{Q\in{\mathcal D}} c_Q \frac{\int_Q \psi(x) \, d\mu(x)}{\mu(Q)} \leq
||g||_{L^{q'}(d\mu)}^{p'}\leq C \, ||\psi||_{L^\frac{q'}{p'}(d\mu)}. $$ Using
duality again we get:  
$${\mathcal W}_{K, \,\sigma}^{\mathcal D} [\mu]=
\sum_{Q\in{\mathcal D} } \frac{c_Q}{\mu(Q)}\chi_Q\in
L^\frac{q(p-1)}{p-q}(d\mu),$$ which  is condition (b).   

Next we show  part (b).   Theorem \ref{theoremdiscrete} and
(\ref{3.2})  give that (\ref{trace}) holds if
  for any $g\in L^{q'}(d\mu)$, $g\geq0$, $$ \int_{\R^n}  {\mathcal W}_{K,\,
\sigma}^{\mathcal D} [gd\mu](x) \, g(x) \, d\mu(x)\simeq \int_{\R^n}
T_{K_{\mathcal D}} [gd\mu]^{p'}(x) \, d\sigma(x)\leq
C \, ||g||_{L^{q'}(d\mu)}^{p'}.$$

Since $(K,\sigma)\in\,$\text{DLBO}, we  have:
 $$
{\mathcal W}_{K, \, \sigma}^{\mathcal D} [gd\mu](x)\leq \left(M_\mu^{H\!L,\,d}[g](x)\right)^{p'-1} 
{\mathcal W}_{K, \, \sigma}^{\mathcal D} [\mu](x).$$    H\"older's
inequality  with exponent $r=\frac{q'}{p'-1}$, gives:
\begin{equation*}\begin{split}\int_{\R^n} {\mathcal
W}_{K, \,\sigma}^{\mathcal D} [gd\mu](x) \, g(x) \,  d\mu(x)\leq & 
C\left(\int_{\R^n}\left(M_\mu^{H\!L, \, {\mathcal D} }[g](x)\right)^{q'}dx\right)^r
\\ & \times \left( \int_{\R^n}\left( g(x) 
{\mathcal W}_{K, \,\sigma}^{\mathcal D} [\mu](x)\right)^{r'}d\mu(x)
\right)^\frac1{r'}. \end{split}\end{equation*}

 Using now
   the boundedness of  $M_\mu^{H\!L,\, {\mathcal D}}$  on 
$L^{q'}(d\mu)$  and H\"older's
inequality with $\lambda=\frac{q'}{r'}>1$ for the last integral, we see that the above
integral is bounded by 
$$C||g||_{L^{q'}(d\mu)}^{p'}\left(
\int_{\R^n}  {\mathcal W}_{K, \, \sigma} ^{\mathcal D}
 [\mu](x)^{r'\lambda'}d\mu(x)\right)^\frac1{r'\lambda'}.$$ Since
$r'\lambda'=\frac{q(p-1)}{p-q}$,  the preceding estimate   gives
(\ref{trace}).\qed 

 \begin{corollary}\label{corollary4.1}
Let $K:{\mathcal D}\rightarrow\R^+$. Let $\mu$ and $\sigma$ 
be  locally finite positive Borel
measures on $\R^n$, $1<q<p<+\infty$. Suppose that 
$(K,\sigma)\in\,${\rm \text{DLBO}}. Then the following statements are
equivalent:

{\rm(a)} There exists $C>0$ such that for any $f\geq 0$,
$$\left(\int_{\R^n}\left(
T_{K_{\mathcal D}}[f](x)\right)^q d\mu(x)\right)^\frac1{q}\leq C \, ||f||_p.$$

 {\rm(b)}
$\, \, \displaystyle{{\mathcal W}_{K, \,\sigma}^{\mathcal D}
[\mu]=\sum_{Q\in{\mathcal D}}\sigma(Q)K(Q)\chi_Q\left(\int_Q
\overline{K}(Q)(y)d\mu(y)\right)^{p'-1} \in L^\frac{q(p-1)}{p-q}(d\mu)}$.\qed 
\end{corollary}

In the general situation without assuming that  $(K,\sigma)\in\,$
\text{DLBO}, we can give   some sufficient conditions in
order that the trace inequality hold. We need to introduce another Wolff-type
potential,  \begin{equation}\label{formula0.2}
\overline{\mathcal{W}}_{K,\,\sigma}^{\mathcal D} [\mu](x)=\sum_{Q\in{\mathcal
D}}\sigma(Q)  \overline{K}(Q)(x)\left(\int_{ Q} \overline{K}(Q)(y)d\mu(y)
\right)^{p'-1} , \end{equation}

  Obviously
${\displaystyle K(Q) \chi_Q(x)\leq \overline{K}(Q)(x)},$ so  that 
for any $x\in\R^n$, we have ${\mathcal{W}}_{K,\,\sigma}^{\mathcal D}
[\mu](x)\leq   \overline{\mathcal{W}}_{K,\,\sigma}^{\mathcal D} [\mu](x)$. We
also observe that when $d\sigma=dx$ is Lebesgue measure on $\R^n$ and
$\overline{K}(Q)=\frac1{r_Q^{n-\alpha}}$ is a discrete Riesz kernel, then   
${\mathcal{W}}_{K,\,dx}^{\mathcal D} [\mu] \simeq  
\overline{\mathcal{W}}_{K,\,dx}^{\mathcal D} [\mu]$.

\begin{theorem}\label{theorem1.2}
Let $1\leq q<p<+\infty$,  and let $\mu$ and $\sigma$ be locally
finite positive Borel measures on $\R^n$.  If $\overline{\mathcal
W}_{K, \,\sigma}^{\mathcal D} [\mu]\in L^{\frac{q(p-1)}{p-q}}(d\mu)$, then
there exists $C>0$ such that for any $f\in L^p(d\sigma)$, $f \ge 0$, 
$$\int_{\R^n}
\left(T_{K_{\mathcal D}}[fd\sigma]\right)^q  \,d\mu \leq C \,
\|f\|_{L^p(d\sigma)}^q.$$      \end{theorem}

 {\bf Proof of Theorem \ref{theorem1.2}:}\par
 
 We can assume that  $q> 1$ since in  
 Corollary \ref{corollary0.1} it was proved that for $q=1$ the condition
${\mathcal W}_{K, \,\sigma}^{\mathcal D} [\mu]\in L^1(d\mu)$ is necessary and
sufficient for the trace inequality to hold, and as we have already observed
${\mathcal W}_{K, \,\sigma}^{\mathcal D} [\mu]\leq \overline{\mathcal W}_{K,
\,\sigma}^{\mathcal D} [\mu]$.  Then  $\overline{\mathcal
W}_{K,\,\sigma}^{\mathcal D} [\mu]\in L^{\frac{q(p-1)}{p-q}}(d\mu)$  is
equivalent by duality to  the fact that there exists $C>0$ such that for any
$g\in L^{\frac{q'}{p'}}(d\mu)$, $g \ge 0$, 
 \begin{equation}\label{qgeq1} \sum_{Q\in{\mathcal D}} \sigma(Q)\left(
\int_Q\overline{K}(Q)(y) \, d\mu(y) \right)^{p'-1}
\int_Q\overline{K}(Q)(y) \, g(y) \, d\mu(y)\leq C \,
\|g\|_{L^{\frac{q'}{p'}}(d\mu)}. \end{equation} 
Next, let $\varphi\in L^{q'}(d\mu)$, $\varphi \ge 0$. Theorem
\ref{theoremdiscrete} and the estimate 
${\mathcal W}_{K, \,\sigma}^{\mathcal D} [\mu]\leq
\overline{\mathcal W}_{K, \,\sigma}^{\mathcal D} [\mu]$ give that
\begin{equation*}\begin{split} \int_{\R^n} \left(T_{K_{\mathcal D}}[\varphi
d\mu](x)\right)^{p'}d\sigma(x)\leq & \, C\int_{\R^n}\overline{\mathcal W}_{K,
\, \sigma}^{\mathcal D} [\varphi d\mu](x) \, \varphi(x) \, d\mu(x)\\ = & \, C
\sum_{Q\in{\mathcal D}}\sigma(Q)\left(
\int_Q\overline{K}(Q)(y)\, \varphi(y) \, d\mu(y) \right)^{p'} .
\end{split}\end{equation*}   Applying  H\"older's
inequality and (\ref{qgeq1}), we obtain:
\begin{equation*}\begin{split}&\int_{\R^n} \left(T_{K_{\mathcal
D}}[\varphi d\mu](x)\right)^{p'} \, d\sigma(x)  \\& \leq \,
C\sum_{Q\in{\mathcal D}}\sigma(Q)
\left(\int_Q\overline{K}(Q)(y) \, d\mu(y)\right)^{\frac{p'}{p}}
\int_Q\overline{K}(Q)(y) \, \varphi(y)^{p'} \, d\mu(y)   \\& \leq C \, 
\|\varphi^{p'}\|_{L^{\frac{q'}{p'}}(d\mu)}=C \,
\|\varphi\|_{L^{q'}(d\mu)}^{p'}. \end{split}\end{equation*} 
Duality again gives that  there exists $C>0$ such
that for any $f\in L^p(d\sigma)$, $f\ge 0$, 
$$\int_{\R^n} (T_{K_{\mathcal
D}}[fd\sigma](x))^q \,d\mu(x)\leq C \, \|f\|_{L^p(d\sigma)}^q.\qed$$

\noindent {\bf Remark.} The condition $\overline{\mathcal
W}_{K, \,\sigma}^{\mathcal D} [\mu]\in L^{\frac{q(p-1)}{p-q}}(d\mu)$ is not
necessary in general, as the following example shows for the case $q=1$ and
$p=2$. Let $k(r)=\frac{1}{r^n \log^\beta (C/r)}$ for $0<r\leq 1$, and $k(r)=0$
for $r>1$, where $1<\beta\leq\frac32$ and $C>0$ is big enough so that $k$ is
nonincreasing (more precisely, we need $C\geq e^\frac{\beta}{n}$). Let $Q_0$ be
the unit cube in $\R^n$,  let $\mu$ be  Lebesgue measure restricted to $Q_0$,
and let  $\sigma$ be Lebesgue measure on $\R^n$.  Then $\sum_Q
k(r_Q)\mu(Q)\chi_Q(x)$ is zero unless $x\in Q_0$, and if $x\in Q_0$,  $$\sum_Q
k(r_Q)\mu(Q)\chi_Q(x)=\sum_{l\geq 0} k(\frac1{2^l}) 2^{-ln} =\sum_{l\geq 0}
\frac1{\log^\beta \frac{C}{2^l}} $$ which is convergent, since $\beta>1$.
Consequently  \begin{equation*}\begin{split} {\mathcal E}_{K_{\mathcal
D}, \, dx}[\mu] & = \int_{\R^n} \left( \sum_Q k(r_Q)\mu(Q)\chi_Q(x) \right)^2
dx\\& = \int_{Q_0} \left( \sum_Q k(r_Q)\mu(Q)\chi_Q(x) \right)^2 dx<+\infty.
\end{split}\end{equation*} Corollary \ref{corollary0.1} gives then that the
trace inequality for $q=1$ and $p=2$ holds. 

On the other hand, if $Q\subset Q_0$,  
$${\overline K}(Q)(x)\simeq \frac{1}{r_Q^n}\int_0^{r_Q} k(t) \, t^{n-1}
\, dt\simeq \frac{1}{r_Q^n\log^{\beta-1}\frac{C}{r_Q} },$$
for any $x\in Q$. Hence for $x\in Q_0$, 
\begin{equation*}\begin{split}&
\overline{\mathcal{W}}_{K,\,\sigma}^{\mathcal D} [\mu](x)=\sum_{Q\in{\mathcal
D}}\sigma(Q)  \overline{K}(Q)(x)\left(\int_{ Q} \overline{K}(Q)(y)d\mu(y)
\right)\geq\\& \sum_{x\in Q\subset
Q_0}r_Q^n\left(\frac{1}{r_Q^n\log^{\beta-1}\frac{C}{r_Q}}
\right)^2\mu(Q)=\\&\sum_{l\geq 0} \frac1{2^{ln}}\left(
\frac{2^{ln}}{\log^{\beta-1}(C2^l)}\right)^2\frac1{2^{ln}}\geq C \sum_{l\geq
1} \frac1{l^{2\beta-2}}=+\infty, \end{split}\end{equation*} since
$\beta\leq\frac32.$ Consequently,
$\overline{\mathcal{W}}_{K,\,\sigma}^{\mathcal D} [\mu](x)\notin
L^1(d\mu)$.\qed

\section{Continuous Wolff-type inequalities and application to continuous
trace inequalities} \label{section02}

One of our main goals in this section is to derive the continuous version  of 
Wolff's inequality   from its discrete counterpart. We start with 
some definitions.

Let $k:(0,+\infty)\rightarrow\R^+$ be a  nonincreasing lower
semicontinuous function, and let $\sigma$ be a positive locally finite Borel
measure on $\R^n$. We set $$
\overline{k}(r)(x)=\frac1{\sigma(B(x,r))}\int_0^r
k(l) \, \sigma(B(x,l)) \, \frac{dl}{l},$$ for $x\in\R^n$, $r>0$.

Our first observation is that if $\sigma$ is a doubling measure 
then $\overline{k}(\cdot)(x)$    satisfies a doubling condition. 
\begin{lemma}\label{doublingoverline}
If $\sigma$ is a  doubling measure 
then  there exists $C>0$ such that, for any $x\in\R^n$, $r>0$,
$$\frac1{C}\,  \overline{k}(2r)(x)\leq \overline{k}(r)(x)\leq
C \, \overline{k}(2r)(x).$$ \end{lemma}  {\bf Proof of Lemma
\ref{doublingoverline}:}\par Since $\sigma$ is doubling, $\sigma(B(x,r))\simeq
\sigma(B(x,2r))$. Then there exists $C>0$ such that for any $x\in\R^n$,
$r>0$, $\overline{k}(r)(x)\leq C\overline{k}(2r)(x)$. On the other hand, the
change of variables $l=2s$, together with the fact that $k$ is nonincreasing, 
yields:
\begin{equation*}\begin{split}\overline{k}(2r)(x) &=\frac1{\sigma(B(x,2r))}
\int_0^{2r}k(l) \, \sigma(B( x,l)) \, \frac{dl}{l}\\& \simeq \frac1{\sigma(B(x,
r))}\int_0^r k(2s) \, \sigma(B(x,2s)) \, \frac{ds}{s}\leq C
\overline{k}(r)(x).\end{split}\end{equation*} \qed

The following lemma gives 
another equivalent reformulation of the function $\overline{k}(r)(x)$.

\begin{lemma}\label{reformulation}
If $\sigma$  is  a
doubling measure then there exists $C>0$ such that, for any
$x\in\R^n$, $r>0$, \begin{equation}\label{refor}\frac1{C} \, \overline{k}(r)(x)\leq
\frac1{\sigma(B(x,r))}\int_{B(x,r)}k(|x-y|) \, d\sigma(y)\leq C \, 
\overline{k}(r)(x).\end{equation}
\end{lemma}
 {\bf Proof of Lemma \ref{reformulation}:}\par
 We begin by proving the second inequality. Since $k$ is nonincreasing
and $\sigma$ is a doubling measure, we have: 
\begin{equation*}\begin{split}  & \int_{B(x,r)}k(|x-y|) \, d\sigma(y)=
\sum_{l\geq 0} \int_{B(x,\frac{r}{2^l})\setminus B(x,\frac{r}{2^{l+1}})}
k(|x-y|) \, d\sigma(y) \\& \leq  \sum_{l\geq0} k(\frac{r}{2^{l+1}}) \left(
\sigma(B(x,\frac{r}{2^l}))-\sigma( B(x,\frac{r}{2^{l+1}}))\right) \leq
\sum_{l\geq0} k(\frac{r}{2^{l+1}})   \sigma(B(x,\frac{r}{2^l})) \\& \leq C
\sum_{l\geq0}\int_{\frac{r}{2^{l+2}}}^{\frac{r}{2^{l+1}}}
k(\frac{r}{2^{l+1}})\sigma(B(x,\frac{r}{2^{l+1}}))\frac{ds}{s} \\& \leq C 
\sum_{l\geq0}\int_{\frac{r}{2^{l+2}}}^{\frac{r}{2^{l+1}}} 
k(s)\sigma(B(x,s))\frac{ds}{s} \leq C\int_0^\frac{r}{2} k(s)
\sigma(B(x,s))\frac{ds}{s}\leq C\overline{k}(r)(x) \, \sigma(B(x,r)).  
\end{split}  \end{equation*}    
To prove the first inequality in (\ref{refor}), we recall that since 
$\sigma$ is a doubling measure, it follows   that there exists
$C>0$ such that for any $x\in\R^n$, $t>0$, $l\in \Z$, 
$$\sigma(B(x,\frac{r}{2^{l}}))\leq C\sigma\left(B(x,\frac{r}{2^{l}})\setminus
B(x,\frac{r}{2^{l+1}})\right) \leq C \sigma(B(x,\frac{r}{2^{l}})).$$ 
Consequently,  \begin{equation*}\begin{split}  &\int_0^r k(s)
\sigma(B(x,s))\frac{ds}{s}=
\sum_{l\geq0}\int_{\frac{r}{2^{l+1}}}^{\frac{r}{2^{l}}}k(s)
\sigma(B(x,s))\frac{ds}{s}\\& \leq \sum_{l\geq0}
\int_{\frac{r}{2^{l+1}}}^{\frac{t}{2^{l}}}k(\frac{r}{2^{l+1}})
\sigma(B(x,\frac{r}{2^{l}}))\frac{ds}{s}\\&  \leq 
C\sum_{l\geq0}k(\frac{r}{2^{l+1}})\sigma\left(B(x,\frac{r}{2^{l}})\setminus
B(x,\frac{r}{2^{l+1}})\right)\leq C
\int_{B(x,\frac{r}2)}k(|x-y|)d\sigma(y).\qed \end{split}\end{equation*} 
   
\begin{lemma}\label{relationship}  Let $k:(0,+\infty)\rightarrow\R^+$ be a
nonincreasing lower semicontinuous function, and let $\sigma$ be a locally
finite positive Borel measure on $\R^n$ satisfying a doubling condition. There
exists $C>0$ such that if  $Q\in{\mathcal D}$ and $x\in Q$,  then 
$$\frac1{C} \, \overline{k}(r_Q)(x)\leq \frac1{\sigma(Q)}\sum_{Q'\subset
Q}k(r_{Q'})\sigma(Q') \chi_{Q'}(x)\leq C \, \overline{k}(r_Q)(x).$$  
In other words, 
if $\overline{K}(Q)(x)$ is the function associated with $K(Q)=k(r_Q)$ and
$\sigma$ is doubling, then $\overline{K}(Q)(x)\simeq \overline{k}(r_Q)(x)$, 
for $x\in Q$. 
\end{lemma} {\bf Proof of Lemma \ref{relationship}:}\par

  Observe  that  $\sigma$ satisfies a doubling
 condition. Hence,  if $Q\in{\mathcal D}$ and $x\in Q$, then
$\sigma(B(x,r_Q))\simeq  \sigma(Q).$  

For any $Q\in{\mathcal D}$,   $x\in Q$ and $l\geq 0$ there exists a unique 
cube $Q_l$ in ${\mathcal D}$ such that $x\in Q_l$, $Q_l\subset Q$ and
$r_{Q_l}=\frac{r_Q}{2^{l}}$. Hence $Q_l\subset B(x, c\frac{r_Q}{2^{l+1}})$,
where $c>0$ is a fixed constant which depends  only on $n$. 

Since $ \sigma(Q_l)\simeq \sigma(B(x, c\frac{r_Q}{2^{l+1}}))$, we have
\begin{equation*}\begin{split}&
\frac1{\sigma(Q)}\sum_{Q'\subset Q}k(r_{Q'}) \sigma(Q')
 \chi_{Q'}(x)=\\&\frac1{\sigma(Q)}\sum_{x\in Q'\subset Q} k(r_{Q'})
\sigma(Q') \simeq \frac1{\sigma(Q)}\sum_{l\geq0}
k(\frac{r_Q}{2^l})\sigma(B(x,c\frac{r_Q}{2^{l}})). \end{split}\end{equation*} 
 But since $\sigma$ satisfies a doubling condition,
$\sigma(B(x,c\frac{r_Q}{2^l}))\simeq \sigma(B(x,\frac{r_Q}{2^l}))$ for any
$l\geq0$. Thus the last sum is bounded above by 
$$C\frac1{\sigma(B(x,r_Q))}\int_0^{r_Q}k(t)\sigma(B(x,t))\frac{dt}{t}=C
\overline{k}(r_Q)(x).$$ On the other hand, $$\frac1{\sigma(Q)}\sum_{l\geq0}
k(\frac{r_Q}{2^l})\sigma(B(x,c\frac{r_Q}{2^{l}}))\succeq
\overline{k}(2r_Q)(x).$$ Since by Lemma \ref{doublingoverline},
$\overline{k}(2r_Q)(x)\simeq \overline{k}(r_Q)(x)$, we obtain the lower 
estimate.\qed

 As for the discrete version, we will restrict ourselves to functions $k$ and
measures $\sigma$ satisfying an extra assumption analogous to property DLBO,
namely that there exists $C>0$ such that for any $x\in\R^n$, $r>0$, and $z\in
B(x,r)$,   \begin{equation}\label{continuousassumption} 
\frac1{C}\overline{k}(r)(x)\leq \overline{k}(r)(z)\leq C \overline{k}(r)(x). 
\end{equation}    In this case we will say that the pair $(k,\sigma)$
satisfies the property of the logarithmic bounded oscillation, or  
simply write $(k,\sigma) \in \text{LBO}$.  The following lemma shows the
relationship between the  LBO and  DLBO properties. The proof is an
immediate consequence of Lemma \ref{relationship}.   
\begin{lemma}\label{continuousdyadic}  Let $k:(0,+\infty)\rightarrow\R^+$ be a
nonincreasing lower semicontinuous function, and let $\sigma$ be a locally
finite positive Borel measure on $\R^n$ satisfying a doubling condition.
Assume that $(k,\sigma)\in${\rm\text{LBO}}. If we set 
$K(Q)=k(r_Q)$, then $(K,\sigma)\in${\rm\text{DLBO}}.\qed  \end{lemma}    
 We observe that the above lemma can be refined in the following sense:  
if for any $z\in\R^n$ and $Q\in{\mathcal D}_z$, where ${\mathcal D}_z$ denotes 
the shifted dyadic lattice ${\mathcal D}_z={\mathcal D}+z$,  we set
$K^z(Q)=k(r_Q)$, then $(K^z,\sigma)\in${\rm\text{DLBO}}, with constants that
do not depend on $z\in\R^n$. We will use this observation later on.     
We first check that the  examples that we  considered in Proposition
\ref{examples}  have continuous analogues satisfying property LBO.    
\begin{proposition}\label{continuousexamples} {\rm(i)} Suppose that
$d\sigma=dx$ is Lebesgue measure on $\R^n$, and suppose that
$k:(0,+\infty)\rightarrow\R^+$ is a nonincreasing lower semicontinuous
function. Then  $(k,dx)\in${\rm\text{LBO}}.        

{\rm(ii)} Suppose that $0<\alpha<n$, and that for any $r>0$,  
$\displaystyle{k(r)=\frac1{r^{n-\alpha}}}$, i.e., $k(|x-y|)$ is the Riesz
kernel on $\R^n$. Suppose     $\sigma$ is a  positive 
doubling  measure on $\R^n$, and   $\sigma \in RD_\gamma$  for some
$\gamma>n-\alpha$, i.e., there exists $C>0$ such that for any $x\in\R^n$,
$A>0$,  $r>0$,  
\begin{equation}\label{reversedoubling}   \sigma(B(x,Ar))\geq
CA^\gamma\sigma(B(x,r)).  \end{equation}  Then 
$(k,\sigma)\in${\rm\text{LBO}}.  \end{proposition}  {\bf Proof of
Proposition \ref{continuousexamples}:}\par   Statement (i) is immediate since
   $\sigma(B(x,t))\simeq t^n$ if $d\sigma=dx$ is  Lebesgue measure on $\R^n$.
Hence  $${\overline{k}}(r)(x)=\frac1{r^n} \int_0^r
k(t) \, t^{n-1}dt$$  is a radial function which obviously
satisfies property LBO.    

Let us show (ii). If
$k(r)=\frac1{r^{n-\alpha}}$ and $\sigma$ is a doubling measure on  $\R^n$
such that $\sigma \in RD_\gamma$,  with
$\gamma>n-\alpha$, then for any $x\in\R^n$, $r>0$, 
\begin{equation*}\begin{split}&  \int_0^r\frac1{t^{n-\alpha}}
\sigma(B(x,t))\frac{dt}{t}=\sum_{l\geq0}\int_{
\frac{r}{2^{l+1}}}^{\frac{r}{2^l}}\frac1{t^{n-\alpha}}\sigma(B(x,t))
\frac{dt}{t}\\& \leq C\sum_{l\geq0}\frac{2^{l(n-\alpha)}}{r^{n-\alpha}}
\sigma(B(x,\frac{r}{2^l})).  \end{split}\end{equation*}  
The fact that
$\sigma$ satisfies a reverse doubling condition gives 
$\sigma(B(x,\frac{r}{2^l}))\leq \frac{C}{2^{l\gamma}}\sigma(B(x,r))$, and
consequently that the above sum is bounded by
$$\left(\sum_{l\geq0}\frac1{2^{l(\gamma-(n-\alpha))}}\right)
\frac{\sigma(B(x,r))}{r^{n-\alpha}}.$$  Since $\gamma-(n-\alpha)>0$, we have
that $\overline{k}(r)(x)\leq C\frac1{r^{n-\alpha}}$.    The fact that $\sigma$
is a doubling measure, gives that 
$$\int_0^r\frac1{t^{n-\alpha}}\sigma(B(x,t))\frac{dt}{t}\geq
\int_{\frac{r}2}^r\frac1{t^{n-\alpha}}\sigma(B(x,t))\frac{dt}{t}\geq
C\frac{\sigma(B(x,r))}{r^{n-\alpha}},$$  and consequently that
$\overline{k}(r)(x)\geq C\frac1{r^{n-\alpha}}$.\qed   

\noindent {\bf Remark.} Observe
that the two examples considered in the above proposition, satisfy a stronger
property, namely that for any $x,z\in\R^n$, $r>0$, $\overline{k}(r)(x)\simeq
\overline{k}(r)(z)$.

 We next define a continuous Wolff-type potential. 
 If $x\in\R^n$, we consider
 \begin{equation}\label{wolffcontinuous}
 {\mathcal W}_{k,\,\sigma}[\mu](x)=\int_0^{+\infty} k(r)\sigma(B(x,r))
 \left(\int_{B(x,r)} \overline{k}(r)(y)d\mu(y)\right)^{p'-1}\frac{dr}{r}.
 \end{equation}
 
 Observe that if $(k,\sigma)\in\text{LBO}$, then the above expression is 
equivalent to
 $$
 {\mathcal W}_{k,\,\sigma}[\mu](x)=\int_0^{+\infty} k(r)\sigma(B(x,r))
\overline{k}(r)(x)^{p'-1} \mu (B(x,r))^{p'-1}\frac{dr}{r}.
 $$
 For technical reasons we will also introduce the truncated Wolff-type
potentials.  If $x\in\R^n$, $R>0$,
 \begin{equation}\label{truncatedwolffcontinuous}
 {\mathcal W}_{k,\,\sigma}^R[\mu](x)=\int_0^{R} k(r)\sigma(B(x,r))
 \left(\int_{B(x,r)} \overline{k}(r)(y)d\mu(y)\right)^{p'-1}\frac{dr}{r}.
 \end{equation}
 If $\mu$ and $\sigma$ are positive locally finite measures on $\R^n$, and
$1<p<+\infty$, the energy associated with $k$ and $\sigma$ is given by 
\begin{equation}\label{generalaenergy}  {\mathcal E}_{k,\,
\sigma}[\mu]=\int_{\R^n}\left( T_k[\mu](x)\right)^{p'}d\sigma(x). 
\end{equation}  
 The following proposition gives a pointwise relationship between the dyadic 
Wolff potential and its continuous version. 
 \begin{proposition}\label{pointwise}
 Let  $k:(0,+\infty)\rightarrow\R^+$
be a  nonincreasing lower semicontinuous function. Let 
$1<p<+\infty$,  and let $\sigma$  be a positive
locally finite Borel measure  on $\R^n$. Suppose also that $\sigma$ satisfies a
doubling condition.  Let $\overline{K}(Q)(x)$ be the function associated with
$K(Q)=k(r_Q)$. Then there exist constants $c, \, C>0$ such that for any
positive locally finite Borel measure $\mu$ on $\R^n$, and $x\in\R^n$,  $$\sum_{Q\in{\mathcal D}}k(cr_Q)\sigma(Q)\chi_Q(x)\left(
\int_Q\overline{ K }(Q)(y)d\mu(y)\right)^{p'-1}\leq C \, {\mathcal
W}_{k,\,\sigma}[\mu](x). $$  \end{proposition}  {\bf Proof of Proposition
\ref{pointwise}:}\par

  If $x\in\R^n$ and $l\in\Z$, there exists a unique cube $Q_l$ in ${\mathcal D}$ such 
that $x\in Q_l$ and $r_{Q_l}=2^{l}$.  Hence $Q_l\subset B(x, \frac{c}2
2^{l})$ where $c>0$ is a fixed constant which  depends only on $n$. 
Applying Lemma \ref{relationship} we
obtain:
   \begin{equation*}\begin{split}&  \sum_{Q\in{\mathcal
D}}k(cr_Q)\sigma(Q)\chi_Q(x)\left( \int_Q\overline{ k
}(r_Q)(y)d\mu(y)\right)^{p'-1}\\&=\sum_{l\in\Z} k(c2^{l}) \sigma(Q_l) \left(
\int_{Q_l} \overline{k}(r_{Q_l})(y)d\mu(y)\right)^{p'-1}\\& \leq C
\sum_{l\in\Z} k(c2^l) \sigma(B(x,\frac{c}2 2^l)) \left( \int_{B(x,\frac{c}2
2^l)} \overline{k}(2^l)(y)d\mu(y)\right)^{p'-1} \\& \leq C \, \sum_{l\in\Z}
\int_{\frac{c}22^l}^{c2^l}k(c2^l)\sigma
(B(x,2^l))\left(\int_{B(x,\frac{c}22^l)}
\overline{k}(2^l)(y)d\mu(y)\right)^{p'-1}\frac{dt}{t}\\& \leq C \,
\sum_{l\in\Z} \int_{\frac{c}2 2^l}^{c2^l}k(t)\sigma(B(x,t))
\left(\int_{B(x,t)}\overline{k}(t)(y)d\mu(y)\right)^{p'-1}\frac{dt}{t}\leq C
{\mathcal W}_{k,\,\sigma}[\mu](x).\qed  \end{split}\end{equation*} 

We now state a continuous version of Wolff's  theorem.

\begin{theorem}\label{theorem2} Let  $k:(0,+\infty)\rightarrow\R^+$
be a  nonincreasing lower semicontinuous function. Let 
$1<p<+\infty$, and  let $\sigma$ be  a positive
locally finite Borel measure on $\R^n$. Suppose that $\sigma$  satisfies a 
doubling condition and that  $(k,\sigma)\in${\rm\text{LBO}}.  Then for any
positive Borel measure $\mu$ on $\R^n$, 
\begin{equation}\label{wolffgeneral}{\mathcal E}_{k, \, \sigma}[\mu]\simeq
\int_{\R^n}{\mathcal W}_{k, \, \sigma}[\mu] \, d\mu,\end{equation} with
constants that may depend on $k$ and $\sigma$, but not on $\mu$. \end{theorem}
{\bf Proof of Theorem \ref{theorem2}:}\par 

For $R>0$, we define the truncated operator $T_{k}^R$ by
$$
T_k^R[\mu](x)= \int_{|x-y|\leq R}k(|x-y|) \, d\mu(y),$$
where $\mu$ is a positive locally finite Borel measure on $\R^n$. The usual 
Fefferman-Stein argument (see \cite{sawyer} and also Lemma 2.2 in \cite{sawyerwheeden}) shows, using 
 the fact that $k$ is nonincreasing, that $T_k^R[\mu] (x)$ is pointwise 
bounded by the average of the shifted dyadic potentials $T_{\tilde K_{{{\mathcal D}_z}}}[\mu](x)$ associated  with  $\tilde k(r) = k(\frac{r}4)$ and $\tilde{K}(Q)=\tilde{k}(r_Q)$, defined by $T_{\tilde K_{{\mathcal D}_z}} [\mu](x)=\sum_{Q\in{\mathcal
D}}k(\frac{r_Q}2) \mu(Q+z)\chi_{Q+z}(x)$. That is, we have that there exists $j_0\in\Z^+$, $C>0$ such that for any $j\in\Z$,
\begin{equation}\label{operatoraverage} T_k^{2^j}[\mu] (x)\leq
\frac{C}{2^{jn}}\int_{|z|\leq 2^{j+j_0}}T_{\tilde K_{{\mathcal D}_z}} [\mu](x)dz.
\end{equation}
Indeed,  fix $j_0$ such that $2^{j_0}>2\sqrt{n}+1$. Then for $x\in
B_j=B(0,2^j)$, and $l\leq j$, we denote by $\Omega$ the set of points $z\in
B_{j+j_0}$, for which there exists $Q\in{\mathcal D}$,
$r_Q=2^{l+1}$, and $B(x,2^j)\subset Q+z$. It is geometrically evident 
that \begin{equation}|\Omega_l|\geq C |B_{j+j_0}|\simeq 2^{jn}.
\label{equation8**}\end{equation}  
Next the fact that $k$ is nonincreasing gives
\begin{equation*}\begin{split}&
T_k^{2^j}[\mu] (x)=\int_{|x-y|\leq 2^j}k(|x-y|)d\mu(y) \\& \leq
 \sum_{l\leq j} k(2^{l-1})\mu(B(x,2^l)).
 \end{split}\end{equation*}
 Applying (\ref{equation8**}) to $l\leq j$ and $x\in B_j$, we obtain:
  \begin{equation*}\begin{split}
   \mu(B(x,2^l))& \leq \frac1{|\Omega_l|} \int_{\Omega_l}\sum_{r_{Q+z}=2^{l+1}}
\mu(Q+z)\chi_{Q+z}(x)dz\\&  \leq  \frac{C}{2^{jn}}
\int_{B_{j+j_0}}\sum_{r_{Q+z}=2^{l+1}} \mu(Q+z)\chi_{Q+z}(x)dz.  
\end{split}\end{equation*}   Altogether, we deduce that
  \begin{equation*}\begin{split}
T_k^{2^j}[\mu] (x)& \leq \frac{C}{2^{jn}} 
\int_{B_{j+j_0}}\sum_{r_{Q+z}=2^{l+1}}k(\frac{r_Q}4) 
\mu(Q+z)\chi_{Q+z}(x)dz \\& \leq \frac{C}{2^{jn}}\int_{|z|\leq
2^{j+j_0}}T_{\tilde K_{{\mathcal D}_z}} [\mu](x)dz, \end{split}\end{equation*}
which  gives (\ref{operatoraverage}).
Now H\"older's inequality together with (\ref{operatoraverage}) gives that
for any $R>0$,  \begin{equation*}\begin{split}
\int_{\R^n}T_k^R[\mu](x)^{p'} \, d\sigma(x) & \leq 
\frac{C}{R^n}\int_{|z|\leq cR}\int_{\R^n}T_{\tilde K_{{\mathcal D}_z}} 
[\mu]^{p'}(x) \, d\sigma(x) \, dz\\& \leq C\sup_{z} \int_{\R^n}T_{\tilde
K_{{\mathcal D}_z}} [\mu]^{p'}(x) \, d\sigma(x). \end{split}\end{equation*}
Applying the dyadic Wolff inequality proved in Theorem \ref{theoremdiscrete} and Lemma \ref{relationship}, 
we obtain that the above expression is bounded by 
$$C\sup_z \sum_{Q\in{\mathcal D}_z} \tilde{k}(r_Q)\sigma(Q)\mu(Q)
\left( \int_Q\overline{{\tilde{k}}}(r_Q)(y)d\mu(y)\right)^{p'-1}  .$$
Next, Lemma \ref{doublingoverline}  gives that
$\overline{\tilde{k}}\simeq  \overline{k}$ so the preceding quantity is bounded
by  \begin{equation}\label{constantc}C\sup_z \sum_{Q\in{\mathcal D}_z}
k(\frac{ r_Q}4) \sigma(Q)\mu(Q)\left( \int_Q\overline{ k
}(r_Q)(y)d\mu(y)\right)^{p'-1}.\end{equation}
 We now have to replace  $k(\frac{r_Q}4)$ in the 
last sum by $k(cr_Q)$,  where $c>0$ is the constant given in Proposition
\ref{pointwise}. This is justified in the following lemma.  
\begin{lemma}\label{lemmaconstant}
Let $k:(0,+\infty)\rightarrow\R^+$ be a  nonincreasing lower 
semicontinuous function. Let  $\sigma$  be a locally finite positive
Borel measure  on $\R^n$,  $1<p<+\infty$ and $0<c<+\infty$. Suppose that $\sigma$ satisfies a
doubling condition and that  $(k,\sigma)\in${\rm\text{LBO}}. Then for any
positive Borel measure $\mu$ on $\R^n$, $$\sum_{Q\in{\mathcal
D}}k(cr_Q)\sigma(Q)\overline{k}(r_Q)^{p'-1}\mu(Q)^{p'}\simeq 
\sum_{Q\in{\mathcal D}}k(r_Q)\sigma(Q)\overline{k}(r_Q)^{p'-1}\mu(Q)^{p'},$$
with constants that do not depend on $\mu$.
\end{lemma}    

{\bf Proof of Lemma \ref{lemmaconstant}:}\par We first observe that since $k$ is nonincreasing 
it follows that if $Q\in{\mathcal D}$ and $2Q$ is the unique cube in ${\mathcal
D}$ such that $Q\subset 2Q$ and $r_{2Q}=2r_Q$, then for any $x\in Q$,
\begin{equation}\label{relation1} \overline{K}(2Q)(x) \simeq
\overline{K}(Q)(x). \end{equation}      Indeed, if $x\in Q$
    $$  \overline{K}(2Q)(x) =
\frac{\sigma(Q)}{\sigma(2Q)}\overline{K}(Q)(x) + k(2r_Q)  \chi_{2Q}(x)\leq
\overline{K}(Q)(x) +k(r_Q)\chi_Q(x)\leq 2\overline{K}(Q)(x).    $$    
   We also observe that since $\sigma$ satisfies a doubling condition, 
   we have that for any $Q\in{\mathcal D}$, $\sigma(Q)\simeq \sigma(2Q)$, and 
consequently for any $Q\in{\mathcal D}$, $x\in Q$,
   \begin{equation}\label{relation2}
   \overline{K}(Q)(x) =\frac{\sigma(2Q)}{\sigma(Q)}
\left( \overline{K}(2Q)(x)-k(2r_Q)\right)  
\leq C \overline{K}(2Q)(x).
   \end{equation}
   
   Next, the fact that $(k,\sigma)\in LBO$ implies by Lemma \ref{continuousdyadic} that $(K,\sigma)\in DLBO$, and consequently that $\overline{K}(2Q) \simeq
\overline{K}(Q)$. Since Lemma \ref{relationship} shows that $\overline{K}(Q)\simeq \overline{k}(r_Q)$, we begin showing that
   \begin{equation}\label{estimate*}\sum_{Q\in{\mathcal
D}}k(cr_Q)\sigma(Q)\overline{K}(Q)^{p'-1}\mu(Q)^{p'} \le C \, 
\sum_{Q\in{\mathcal D}}k(r_Q)\sigma(Q)\overline{K}(Q)^{p'-1}\mu(Q)^{p'}.\end{equation}
   
   Since $k$ is a nonincreasing function, we can assume 
without loss of generality that 
$c=\frac1{2^l}$,  $l\geq 0$.  We  have
   \begin{equation}\label{constant1}
   \sum_{Q\in{\mathcal D}} \sigma(Q)k(\frac1{2^l}r_Q) 
\overline{K}(Q)^{p'-1}\mu(Q)^{p'} \simeq
   \sum_{Q\in{\mathcal D}} \sigma(2^l Q)k( r_Q) 
\overline{K}(2^lQ)^{p'-1}\mu(2^lQ)^{p'}.
   \end{equation}
   
   But $\mu(2^lQ)= \sum_{Q'\subset 2^l Q,\, r_{Q'}=r_Q}\mu(Q')$, where 
the sum is taken over
all cubes $Q'$ in  ${\mathcal D}$   that are contained in $2^lQ$ and such
that $r_{Q'}=r_Q$. The doubling condition imposed on $\sigma$  gives
$\sigma(2^lQ)\simeq \sigma(Q)\simeq \sigma(Q')$, and (\ref{relation1}) implies 
$\overline{K}(2^lQ)\simeq \overline{K}(Q)\simeq\overline{K}(Q')$. Thus
the left-hand side of (\ref{constant1}) is bounded above by
$$C\sum_{Q\in{\mathcal D}} \sigma(Q) k(r_Q)
\overline{K}(Q)^{p'-1}\mu(Q)^{p'},$$      
and we have (\ref{estimate*})
for $c=\frac 1{2^l}$, $l \ge 0$. 
The converse estimate is 
obvious because $k$ is a nonincreasing function.\qed

We now complete the proof of the theorem.  Lemma \ref{lemmaconstant} implies 
that we can replace $\frac14$ in (\ref{constantc}) by any positive constant,
e.g., by the constant $c>0$ given in Proposition \ref{pointwise}. Therefore,
(\ref{constantc}) is bounded by $$  C \, \sup_z \sum_{Q\in{\mathcal D}_z} k(
cr_Q )\sigma(Q) \overline{k} (r_Q)^{p'-1}   \mu(Q)^{p'}. $$ Finally, the
pointwise inequality obtained in Proposition \ref{pointwise}  gives that the
above expression is bounded by $C \, \int_{\R^n}{\mathcal
W}_{k,\,\sigma}[\mu](x)d\mu(x),$ which shows that
$$ {\mathcal E}_{k, \,\sigma}[\mu] \leq C\int_{\R^n}{\mathcal
W}_{k,\,\sigma}[\mu](x)d\mu(x).$$ The converse estimate is proved in 
the following lemma.

\begin{lemma} \label{lemma5}  Let  $k:(0,+\infty)\rightarrow\R^+$
be a  nonincreasing lower semicontinuous function. Let  
$1<p<+\infty$,  and let  $\sigma$ be positive
locally finite Borel measure on $\R^n$. Suppose that
$(k,\sigma)\in \,${\rm\text{LBO}}. Then  there exists $C>0$ such that  for any
positive Borel measure $\mu$ on $\R^n$,  
\begin{equation}
{\mathcal E}_{k,\,\sigma}[\mu]\geq C \int_{\R^n}
 {\mathcal W}_{k,\,\sigma}[\mu](x) \, d\mu(x).\label{equation10}
\end{equation}
\end{lemma}
{\bf Proof of Lemma \ref{lemma5}:}\par Let $M^{H\!L}_\sigma$ denote the
centered  Hardy-Littlewood maximal  function with respect to the measure
$\sigma$. We  have \begin{equation}\label{fractionalmaximal}  M^{H\!L}_\sigma
[T_k[\mu]](x)\geq C \, M_k[\mu](x),\end{equation} where $M_k$ is the following 
maximal function associated with the kernel $\bar k$:
$$M_k[\mu](x)=\sup_{r>0} \, \overline{k}(r)(x) \, \mu(B(x,r)).$$ Indeed,
Fubini's theorem  gives 
\begin{equation*}\begin{split} & M^{H\!L}_\sigma 
[T_k[\mu]](x)\geq \frac1{\sigma(B(x,2r))}\int_{B(x,2r)}T_k[\mu](y) \,
d\sigma(y) \\& \geq  \frac{C}{\sigma(B(x,r))}\int_{B(x,2r)}\left(
\int_{B(x,r)}k(|y-z|) \, d\mu(z)\right) \, d\sigma(y)\\& \geq 
\frac{C}{\sigma(B(x,r))}\int_{B(z,r)}\left(
\int_{B(x,r)}k(|y-z|) \, d\sigma(y)\right) \, d\mu(z)\\ &\geq    
\frac{C}{\sigma(B(x,r))}\int_{B(x,r)} \overline{k}(r)(z) \, \sigma(B(z,r))
\, d\mu(z). \end{split}\end{equation*}  
Since  $(k,\sigma)\in
\,${\rm\text{LBO}}, we can replace  ${\overline{k}}(r)(z)$  in the last
integral by ${\overline{k}}(r)(x)$. Next, the fact that $\sigma$ is a doubling
measure   gives that  the above integral is  bounded from below by
${\overline{k}}(r)(x) \, \mu(B(x,r))$, and consequently that
$M^{H\!L}_\sigma [T_k[\mu]](x)\geq C \, M_k[\mu](x)$.

Since $M^{H\!L}_\sigma$ is a bounded operator on $L^{p'}(d\sigma)$, and 
$M^{H\!L}_\sigma [T_k[\mu]]\geq T_k[\mu]$ a.e. (with respect to $d\sigma$), Fubini's
theorem gives
\begin{equation*}\begin{split} ||T_k[\mu]||_{L^{p'}(d\sigma)}^{p'}\geq \, & C
\, || M^{H\!L}_\sigma  [T_k[\mu]]||_{L^{p'}(d\sigma)}^{p'} \\ = \, &
C \, \int_{\R^n}  M^{H\!L}_\sigma   [T_k[\mu]] (x) \left(M^{H\!L}_\sigma [
T_k[\mu]]\right)^{p'-1}(x) \, d\sigma(x) \\  \geq \, &
C\, \int_{\R^n} T_k[\left(M^{H\!L}_\sigma [ T_k[\mu]]\right)^{p'-1}](y) \, 
d\mu(y). \end{split}\end{equation*} But 
\begin{equation*}\begin{split}
&T_k[\left(M^{H\!L}_\sigma [T_k[\mu]]\right)^{p'-1}](y)=\int_{\R^n} k(|x-y|) 
\left( M^{H\!L}_\sigma  [T_k\mu]]\right)^{p'-1}(x) \, d\sigma(x)  \\&=
\sum_{l\in\Z} \int_{\frac1{2^{l+1}}\leq|x-y|<\frac1{2^l}}k(|x-y|) 
\left( M^{H\!L}_\sigma  [T_k\mu]]\right)^{p'-1}(x) \, d\sigma(x)\\& \geq
\sum_{l\in\Z}\int_{\frac1{2^{l+1}}\leq|x-y|<\frac1{2^l}} k(\frac1{2^{l}})
\left( M^{H\!L}_\sigma [T_k\mu]]\right)^{p'-1}(x) \, d\sigma(x). 
\end{split}\end{equation*}  Next
  (\ref{fractionalmaximal}) shows that $M^{H\!L}_\sigma  [T_k\mu]](x)\geq C
\overline{k}(\frac1{2^{l-2}})(x) \mu(B(x,\frac1{2^{l-2}}))$. Since  for any
$y\in\R^n$ such that $\frac1{2^{l+1}}\leq |x-y|<\frac1{2^l}$ we have 
$B(y,\frac1{2^{l-1}})\subset B(x,\frac1{2^{l-2}})$, and by 
(\ref{continuousassumption}) $\overline{k}(\frac1{2^{l-2}})(x)\simeq
\overline{k}(\frac1{2^{l-2}})(y)$  
it follows that the above sum 
 is bounded from below by 
$$ C \, \sum_{l\in\Z}
k(\frac1{2^{l}}) \, \overline{k}(\frac1{2^{l-2}})(y)^{p'-1} \, 
\mu(B(y,\frac1{2^{l-1}}))^{p'-1} \,  \sigma( \{x\in\R^n\,;\,
\frac1{2^{l+1}}\leq |x-y|<\frac1{2^l}\}).$$  The fact that $\sigma$ satisfies
a doubling condition gives, as we have already pointed out earlier, that
$\sigma( \{x\in\R^n\,;\, \frac1{2^{l+1}}\leq |x-y|<\frac1{2^l}\})\simeq
\sigma(B(y,\frac1{2^l}))$. Altogether we have  that the above sum is in its
turn bounded from below by \begin{equation*}\begin{split}&C\sum_{l\in\Z}
k(\frac1{2^{l}}) \, \overline{k}(\frac1{2^{l-2}})(y)^{p'-1} \, 
\mu(B(y,\frac1{2^{l-1}}))^{ p'-1} \, \sigma(B(y,\frac1{2^l})) \\& \geq
C \, \sum_{l\in\Z}\int_{\frac1{2^l}}^{\frac1{2^{l-1}}}
k(\frac1{2^l}) \, \overline{k}( \frac1{2^{l-2}})(y)^{p'-1} \,  \mu(B(y,
\frac1{2^{l-1}}))^{p'-1} \, \sigma(B(y,\frac1{2^l})) \, \frac{dt}{t} \\& \geq
\, C \sum_{l\in\Z}\int_{\frac1{2^l}}^{\frac1{2^{l-1}}} k(t) \, \overline{k}(
t)(y)^{p'-1} \, \mu(B(y, t))^{p'-1} \, \sigma(B(y,t)) \, \frac{dt}{t} 
\end{split}\end{equation*} 
Thus $$ ||T_k[\mu]||_{L^{p'}(d\sigma)}^{p'}\geq
C\int_{\R^n}\int_0^{+\infty} k(t) \, \overline{k}(t)(y)^{p'-1} \, 
\mu(B(y,t))^{p'-1} \, \sigma(B(y,t)) \, \frac{dt}{t} \, d\mu(y).\qed$$

\noindent{\bf Remark.} As in the dyadic case, Theorem \ref{theorem2} 
is no longer true with the nonlinear potential 
$\overline{{\mathcal W}}_{k, \, \sigma}[\mu]$ in place of 
${\mathcal W}_{k, \, \sigma}[\mu]$ even in the case 
$d \sigma = dx$ when
$$\overline{{\mathcal W}}_{k, \, dx}[\mu] (x) = \int_0^{+\infty}  \, 
r^n \, 
\overline{k}(r)^{p'} \, \mu (B(x,r))^{p'-1} \, \frac{dr}{r}$$
and
$\overline{k}(r) = r^{-n} \int_0^r k(s) \, s^{n-1} \, ds.$

Indeed, let $q=1$ and  $p=2$, and let  $d \mu = \chi_{B_0} \, dx$ where
$B_0=B(0,1)$ is the unit  ball in $\R^n$. Then for $k(r) = 
r^{-n} \,
\log^{-\beta} (C/r)$ if $0<r<1$, and $k(r)=0$ if $r\ge 1$, 
  where  $1<\beta\le
\frac 3 2$ and $C\ge e^{\beta/n}$, one has  as in  the example 
at the
end of Sec.~2  that ${\mathcal W}_{k, \, dx}[\mu]$ is uniformly bounded 
and hence  by Theorem \ref{theorem2}, ${\mathcal E}_{k,\, dx}[\mu]<
+\infty$. On the other hand,  $\overline{{\mathcal W}}_{k, \,
dx}[\mu]\equiv +\infty$ on $B_0$ and so 
$\int_{\R^n} \overline{{\mathcal W}}_{k, \, dx}[\mu] 
\, d \mu=+\infty.$
\qed
\bigskip

As in the discrete case, the continuous Wolff-type theorem that we 
have just proved  yields a characterization of the corresponding trace
inequality.  \begin{corollary}\label{continuousq=1}
 Let  $k:(0,+\infty)\rightarrow\R^+$
be a  nonincreasing lower semicontinuous function. Let 
$1<p<+\infty$,  and let $\mu$ and $\sigma$ be positive
locally finite Borel measures on $\R^n$. Suppose also that $\sigma$ satisfies
a  doubling condition and that  $(k,\sigma)\in \,${\rm\text{LBO}}. 
Then the following assertions are equivalent: \begin{itemize}
\item[(i)] There exists $C>0$ such that for any $f\in L^p(d\sigma)$, $f \ge
0$,   $$\int_{\R^n}
T_k[f] \, d\mu \leq C \, ||f||_p.$$
\item[(ii)] If ${\mathcal W}_{K,\,\sigma}^{{\mathcal D}_z} [\mu]$ are the
dyadic shifted potentials defined in {\rm(\ref{formula0.1.0})}, 
$$\sup_{z\in\R^n}\int_{\R^n} {\mathcal W}_{K\,\sigma}^{{\mathcal D}_z}[\mu] \,
d\mu<+\infty.$$

 \item[(iii)]
${\mathcal W}_{k,\,\sigma}[\mu] \in L^1(d\mu)$.

\end{itemize}
\end{corollary}
{\bf Proof of Corollary \ref{continuousq=1}:}\par
Using Fubini's theorem and duality it is easy to see  that (i) is equivalent to
the fact that  $||T_k[\mu]||_{L^{p'}(d\sigma)}<+\infty$. By Theorem
\ref{theorem2} this  is in  turn equivalent to ${\mathcal W}_{k,\,\sigma}
[\mu] \in L^1(d\mu)$ which gives the equivalence of (i) and (iii).

Next, Proposition \ref{pointwise}, together with Lemma \ref{continuousdyadic},
  shows that there exists $c>0$ such that for any $z\in\R^n$, $$\sum_{Q\in {\mathcal D}_z} k( cr_Q
)\sigma(Q)  \overline{K} (Q)^{p'-1}   \mu(Q)^{p' } \leq C \int_{\R^n}{\mathcal
W}_{k,\,\sigma}[\mu](x)d\mu(x).$$  
Lemma \ref{lemmaconstant} gives that the constant $c$ in the above  sum 
on the left-hand side 
can be dropped, and consequently, we have that  (iii) $\Rightarrow$ (ii).

For the last implication, (ii)$\Rightarrow$(i), we proceed as in the
proof of Theorem \ref{theorem2}, using the fact that  the truncated operator
$T_{K}^R$   is pointwise 
bounded by the average of the shifted dyadic potentials 
$T_{\tilde K_{{\mathcal D}_z}} [\mu](x)$. H\"older's inequality gives then 
that for any $R>0$, 
$$
\int_{\R^n}T_k^R[\mu](x)^{p'} \, d\sigma(x)\leq C\sup_{z}
 \int_{\R^n}T_{\tilde K_{{\mathcal D}_z}} [\mu](x)^{p'} \, d\sigma(x).
$$
Applying the dyadic Wolff inequality established in Theorem
\ref{theoremdiscrete}  we obtain that the above expression can be bounded by 
$$C\sup_z \sum_{Q\in{\mathcal D}_z} \tilde{k}(r_Q)\sigma(Q)\mu(Q)
\left( \int_Q\overline{{\tilde{k}}}(r_Q)(y)d\mu(y)\right)^{p'-1}  .$$
 Since $\overline{k}(\cdot)$ satisfies a doubling condition, and by Lemma
\ref{lemmaconstant} we can replace $\tilde{k}(r_Q)$ by $k(r_Q)$, the fact that
 (ii) holds gives that $\int_{\R^n}T_k^R[\mu](x)^{p'}d\sigma(x)\leq C$ for any
$R>0$. Lebesgue's monotone convergence theorem  finally gives (i).\qed   

We now consider the trace inequality for $q\geq 1$.
\begin{theorem}\label{continuousqgeq1}
 Let  $k:(0,+\infty)\rightarrow\R^+$
be a  nonincreasing lower semicontinuous function, $1\leq q<p<+\infty$, and 
$\mu$ and $\sigma$  be positive
locally finite Borel measures on $\R^n$. Suppose also that $\sigma$ satisfies
a  doubling condition and that  $(k,\sigma)\in\,${\rm\text{LBO}}. 
Then the following assertions are equivalent: \begin{itemize}
\item[(i)] There exists $C>0$ such that for any $f\in L^p(d\sigma)$, $f\ge 0$, 
 $$\left(\int_{\R^n}
T_k[f]^q \, d\mu \right)^{\frac 1 q} \leq C \, ||f||_p.$$
\item[(ii)] $\sup_{z\in\R^n}\int_{\R^n} {\mathcal W}_{K,\,\sigma}^{{\mathcal D}_z}
[\mu]^\frac{q(p-1)}{p-q} \, d\mu<+\infty$.

 \item[(iii)]
$\int_{\R^n}  {\mathcal W}_{k,\,\sigma}[\mu]^\frac{q(p-1)}{p-q}
 \, d\mu<+\infty$.
\end{itemize}
\end{theorem}
{\bf Proof of Theorem \ref{continuousqgeq1}:}\par
By the last corollary, we may assume that $q>1$. We begin by showing that (i) 
implies (ii). We observe that Theorem \ref{theorem2} applied to the positive
measure $gd\mu$ gives that 
$$||T_k[gd\mu]||_{L^{p'}(d\sigma)}^{p'}\simeq
\int_{\R^n}{\mathcal W}_{K,\,\sigma}[gd\mu](x)g(x)d\mu(x),$$ 
with constants independent of $g$ and $\mu$. But
Proposition \ref{pointwise} together with Lemma \ref{lemmaconstant} 
applied to the shifted lattice ${\mathcal D}_z$ show that the
above integral is bounded from below by 
$$C\sup_{z\in\R^n}\sum_{Q\in {\mathcal
D}_z} k(r_Q)\sigma(Q) \overline{K}(Q)^{p'-1} \left(\int_Q
g(x)d\mu(x)\right)^{p'}.$$

Now we can proceed as in the proof of (a)$\Rightarrow$(b) in Theorem
\ref{theorem4.2} and obtain that (i)$\Rightarrow$(ii).

 Before we present the rest of the proof, we need an estimate similar to  
(\ref{operatoraverage})
for the average of ${\mathcal W}_{K,\,\sigma}^{\mathcal D}$ over the shifts
of the dyadic lattice ${\mathcal D}$. Recall that, for $R>0$, $x\in\R^n$, we
have
 $${\mathcal
W}_{k,\,\sigma}^R[\mu](x)=\int_0^R k(r) \, \sigma(B(x,r))\left(
\int_{B(x,r)}\overline{k}(r)(y) \, d\mu(y)\right)^{p'-1} \, \frac{dr}{r}.$$ 
Since $(k,\sigma)\in\,${\rm\text{LBO}} it follows that the truncated Wolff type
potential can be rewritten in the equivalent form: 
$${\mathcal
W}_{k,\,\sigma}^R[\mu](x)=\int_0^R k(r) \, \sigma(B(x,r))  \, 
\overline{k}(r)^{p'-1}(x) \, \mu( B(x,r))^{p'-1} \, \frac{dr}{r}.$$ We then
have that there exists $j_0\in \Z^+$ and $C>0$ such that for any $j\in\Z$, 
$x\in B_j=B(0,2^j)$,  
\begin{equation}\label{average}{\mathcal
W}_{k,\,\sigma}^{2^j}[\mu](x)\leq \frac{C}{|B_{j+j_0}|}\int_{B_{j+j_0}}
{\mathcal W}_{{\widetilde K},\,\sigma}^{{\mathcal D}_z}[\mu](x)dz,
\end{equation} where ${\widetilde k}(t)=k(\frac{t}4)$. The proof follows that
of (\ref{operatoraverage}) in Theorem \ref{theorem2}. With the  notations used
there,  fix $j_0$ such that $2^{j_0}>2\sqrt{n}+1$. Then for $j\in\Z$, $x\in
B_j=B(0,2^j)$,  and $l\leq j$,  $\Omega_l$ is the set of points $z\in
B_{j+j_0}$, for which there exists $Q\in{\mathcal D}$, $r_Q=2^{l+1}$, and
$I=B(x,2^l)\subset Q+z$. We recall that by (\ref{equation8**}),
$|\Omega_l|\simeq |B_{j+j_0}|\simeq 2^{jn}$.  The fact that
$\overline{k}(\cdot)$ satisfies a doubling condition gives
\begin{equation*}\begin{split}& {\mathcal W}_{k, \, \sigma}^{2^j}
[\mu](x)=\int_0^{2^j} k(r) \, \sigma(B(x,r)) \left(
\int_{B(x,r)} \overline{k}(r)(y) \, d\mu(y)\right)^{p'-1} \frac{dr}{r}\\ & =
\sum_{l\leq j}\int_{2^{l-1}}^{2^l} k(r) \, \sigma(B(x,r))\left(
\int_{B(x,r)}\overline{k}(r)(y) \,  d\mu(y)\right)^{p'-1} \, \frac{dr}{r} \\&
\leq C\sum_{l\leq j}\sigma(B(x,2^l)) \, 
k(2^{l-1}) \, \overline{k}(2^l)(x)^{p'-1} \, \mu(B(x,2^l))^{p'-1}. \end{split}
\end{equation*} Applying (\ref{equation8**}) to $l\leq j$ and $x\in B_j$, we
conclude  that 
\begin{equation*}\begin{split}&\mu(B(x,2^l))^{p'-1}\leq
\frac{1}{|\Omega_l|}\int_{\Omega_l}\sum_{r_{Q+z}=2^{l+1}, \, B(x,2^l)\subset Q+z}
\mu(Q+z)^{p'-1}\chi_{Q+z}(x)dz\\
&\leq \frac{C}{|B_{j+j_0}|}\int_{B_{j+_0}}\sum_{r_{Q+z}=2^{l+1}, \, B(x,2^l)\subset Q+z}
\mu(Q+z)^{p'-1}\chi_{Q+z}(x)dz. \end{split} \end{equation*}
Hence, if $x\in B_j$,
\begin{equation*}\begin{split}& {\mathcal W}_{K, \, \sigma}^{2^j}
[\mu] (x)\leq \frac{C}{|B_{j+j_0}|}\int_{B_{j+j_0}} \sum_{l\leq j}
\sum_{r_{Q+z}=2^{l+2}}\sigma(Q+z)k(\frac{r_Q}4) \overline{k}(r_Q)(x)^{p'-1}\\ 
& \times 
\mu(Q+z)^{p'-1}\chi_{Q+z}(x)dz\leq 
 \frac{C}{|B_{j+j_0}|}\int_{B_{j+j_0}}{\mathcal W}_{\tilde{K}, \, \sigma}
^{{\mathcal D}_z} [\mu] (x) \, dz. \end{split}
\end{equation*}
In the last inequality we have used the estimate ${\overline{\widetilde
k}}(\cdot) \simeq \overline{k} (\cdot)$ which follows from the fact that 
$\overline{k}(\cdot)$  satisfies a doubling condition.

Now we can complete the proof of (ii)$\Rightarrow$(i). 
 Duality and Theorem \ref{theorem2} gives that (i) holds if
  for any $g\in L^{q'}(d\mu)$, $g\geq0$, $$ \int_{\R^n} {\mathcal W}_{k,\,
\sigma}[gd\mu](x) \, g(x) \, d\mu(x)\simeq \int_{\R^n}
T_k [gd\mu]^{p'}(x) \, d\sigma(x)\leq C \, ||g||_{L^{q'}(d\mu)}^{p'}.$$

Now, we consider the translated dyadic Hardy-Littlewood maximal 
function with respect
to $\mu$ given by $$ M_\mu^{H\!L,\,{\mathcal D}_z}h(x)=\sup_{x\in Q+z\,Q\in
{\mathcal D} }\frac1{\mu(Q+z)}\int_{Q+z} |h(y)| \, d\mu(y). $$ We have
that $$ {\mathcal W}_{{\widetilde K}, \, \sigma}^{{\mathcal D}_z}[gd\mu](x)\leq M_\mu^{H\!L,\,{\mathcal D}_z}g(x)^{p'-1}{\mathcal W}_{{\widetilde K},
\, \sigma}^{{\mathcal D}_z}[\mu](x).$$   H\"older's inequality  with
exponent $r=\frac{q'}{p'-1}$, gives
\begin{equation*}\begin{split}&\int_{\R^n} {\mathcal
W}_{{\widetilde K},\, \sigma}^{{\mathcal D}_z}[gd\mu](x)\, g(x) \, d\mu(x)
\\& \leq C \left(\int_{\R^n}M_\mu^{H\!L, {\mathcal D}_z}g(x)^{q'}d\mu
(x)\right)^\frac1{r} \left( \int_{\R^n}\left(  {\mathcal
W}_{{\widetilde K},\,\sigma}^{{\mathcal D}_z}[\mu](x)g(x)\right)^{r'}
d\mu(x)  \right)^\frac1{r'}.\end{split}\end{equation*}

 Using the fact 
   that $M_\mu^{H\!L,{\mathcal D}_z}$ is a bounded operator  on 
$L^{q'}(d\mu)$, and H\"older's
inequality with $\lambda=\frac{q'}{r'}>1$, we have that the above
integral is bounded by $$C||g||_{L^{q'}(d\mu)}^{p'}\left(
\int_{\R^n} {\mathcal W}_{{\widetilde K}, \, \sigma} ^{{\mathcal D}_z 
}[\mu](x)^{r'\lambda'}d\mu(x)\right)^\frac1{r'\lambda'}.$$ 

Since
$r'\lambda'=\frac{q(p-1)}{p-q}$, and we are assuming that (ii)
holds, the last estimate and (\ref{average}) easily  give that
\begin{equation}\label{truncated}\int_{B_j}{\mathcal W}_{k, \, \sigma}^{2^
j}[gd\mu](x) \, g(x) \, d\mu(x)\leq \, C \,
||g||_{L^{q'}(d\mu)}^{p'}\end{equation} once we show that in the expression
\begin{equation}\label{constantcgeneral} \int_{\R^n} {\mathcal
W}_{\widetilde{K}, \,  \sigma}^{{\mathcal D}_z}[\mu](x)^{r'\lambda'} \,
d\mu(x), \end{equation} 
we can replace ${\widetilde{K}}$ by $K$.
This is proved in the following lemma. 
\begin{lemma}\label{lemmaconstantrgeq1}
Let $k:(0,+\infty)\rightarrow\R^+$ be a  nonincreasing 
lower semicontinuous function. 
Let $\sigma$  be locally finite positive Borel measure on $\R^n$,
and let $1<p,r<+\infty$. Assume that $\sigma$ satisfies a doubling condition
and that $(k,\sigma)\in\,${\rm\text{LBO}}. Then  for
any $c>0$ there exists $C>0$ such that  for any positive Borel measure $\mu$ 
on $\R^n$,
\begin{equation*}\begin{split}&\frac1{C}\int_{\R^n}\left(\sum_{Q\in{\mathcal
D}} k(cr_Q)
\, \sigma(Q) \, \overline{K}(Q)^{p'-1}
\, \mu(Q)^{p'-1} \, \chi_Q(x)\right)^r \, d\mu(x)  \\& \leq
\int_{\R^n}\left(\sum_{Q\in{\mathcal
D}} k(r_Q)
\, \sigma(Q) \, \overline{K}(Q)^{p'-1}
\, \mu(Q)^{p'-1} \, \chi_Q(x)\right)^r \, d\mu(x)  \\& \leq C
\int_{\R^n}\left(\sum_{Q\in{\mathcal
D}}
k(cr_Q)
\, \sigma(Q) \, \overline{K}(Q)^{p'-1}
\, \mu(Q)^{p'-1} \, \chi_Q(x)\right)^r \, d\mu(x). \end{split}\end{equation*}

\end{lemma}
   {\bf Proof of Lemma \ref{lemmaconstantrgeq1}:}\par
Note that in the case $r=1$  this lemma coincides with Lemma
\ref{lemmaconstant}.  Since $k$ is nonincreasing,
we can assume without loss of generality that $c=\frac1{2^l}$, $l\geq0$,  and
write $k_l(r)(x)=k(\frac1{2^l}r)(x)$, and $K_l(Q)(x)=k_l(r_Q)(x)$. 
The upper estimate is obvious because $l\geq0$ and $k$ is nonincreasing. The
lower estimate  can be restated  equivalently in the form
\begin{equation}\label{eqv}
\int_{\R^n} ( {\mathcal W}_{K_l, \, \sigma}^{\mathcal D} [\mu](x))^r \,
d\mu(x) \leq \, C \int_{\R^n} ({  {\mathcal W}}_{K, \, \sigma}^{\mathcal D}
[\mu](x))^r \, d\mu(x). \end{equation}

 We have that if  
 $||g||_{L^{r'}(\mu)}\leq 1$, then 
\begin{equation*}\begin{split}&\left[\int_{\R^n} ( {\mathcal W}_{K_l, \,
\sigma}^{\mathcal D} [\mu](x))^r \, d\mu(x)\right]^\frac1{r} \leq \int_{\R^n} 
 {\mathcal W}_{K_l, \,
\sigma}^{\mathcal D} [\mu](x) \, g(x) \, d\mu(x) \\
& \leq C\sum_Q k(\frac1{2^l} r_Q) \, \overline{ k} (r_Q)^{p'-1} \, \sigma(Q)
\,  \mu(Q)^{p'}  \frac 1 {\mu(Q)} \int_Q g(x) \, d
\mu(x).\end{split}\end{equation*}  Similarly to the argument in  Lemma
\ref{lemmaconstant}, we  estimate $\mu(Q)^{p'}$ by     
$C \sum_{Q'}
\mu(Q')^{p'}$ where the sum is taken over the dyadic cubes $Q'$ contained 
in $Q$ such
that $r_{Q'} = 2^{-l} r_Q$  (there are $2^{n l}$ such  $Q'$). 
The doubling condition on 
$\sigma$  gives that $\sigma(Q')\simeq \sigma(Q)$ and
$\overline{k}(r_{Q'})(x)\simeq \overline{k}(r_Q)(x)$. We finally obtain 
$$\left[\int_{\R^n} (
{\mathcal W}_{K_l, \, \sigma}^{\mathcal D} [\mu](x))^r \, d\mu(x)
\right]^\frac1{r}  \leq C \,  \int_{\R^n} {\mathcal W}_{K,\, 
\sigma}^{\mathcal D}
[\mu](x) \,  M^{\mathcal D}_\mu [g](x)  \, d\mu(x).$$ Applying  H\"older's
inequality and the maximal inequality, we get (\ref{eqv}). \qed   

To complete the proof of  (ii)$\Rightarrow$(i)  we let  
$j\rightarrow +\infty$ on the left-hand side of
(\ref{truncated}) and 
apply Lebesgue's monotone convergence theorem.

  The implication (iii)$\Rightarrow$(ii) follows from the pointwise
  estimate given by Proposition \ref{pointwise}, using again the fact that
by Lemma \ref{lemmaconstantrgeq1} we can replace $k(r_Q)$ by $k(c r_Q)$, for
any $c>0$, in the expression $$\int_{\R^n} ({{\mathcal W}}_{K, \,
\sigma}^{\mathcal D} [\mu](x))^r \, d\mu(x).$$

It remains to prove that (ii)$\Rightarrow$(iii). H\"older's inequality
with exponent $\frac{q(p-1)}{p-q}>1$, together with (\ref{average}), 
gives  that \begin{equation*}\begin{split} &\int_{B_k}\left({\mathcal
W}_{k,\,\sigma}^{2^j}[\mu](x)\right)^{\frac{q(p-1)}{p-q}}d\mu(x)\\ & \leq
\frac{C}{|B_{j+3}|}\int_{B_j}\int_{B_{j+3}}\left({\mathcal W}_{\tilde{K},\,
\sigma} ^{{\mathcal D}_z} [\mu](x)\right)^{\frac{q(p-1)}{p-q}}
\, dz \, d\mu(x)\\& =\frac{C}{|B_{j+3}|}\int_{B_{j+3}}\int_{B_j}\left({\mathcal
W}_{\tilde{K},\, \sigma} ^{{\mathcal D}_z
}[\mu](x)\right)^{\frac{q(p-1)}{p-q}}d\mu(x) \, dz.
\end{split}\end{equation*} Now, Lemma \ref{lemmaconstantrgeq1} and  (ii)
easily give that $$\int_{B_k}\left( {\mathcal
W}_{K,\,\sigma}^{2^j}[\mu](x)\right)^{\frac{q(p-1)}{p-q}} \, d\mu(x)\leq C,$$
and letting $j\rightarrow +\infty$, we obtain (iii).\qed


\begin{thebibliography}{CaOrVeW}

\bibitem[{\bf Ad}]{adams} D.R. Adams,   {\it Traces of potentials
arising from translation invariant operators}, Ann. Scuola Norm.
Sup. Pisa Cl. Sci.  {\bf 25} (1971), 203--217.



\bibitem[{\bf AdHe}]{adamshedberg} D.R. Adams and L.I. Hedberg, 
{\it Function
Spaces and Potential Theory}, Springer-Verlag
Berlin--Heidelberg--New York, 1996.

\bibitem[{\bf AiEs}]{aikawa} H. Aikawa and  M. Ess\'en,   {\it Potential
Theory--Selected Topics}, Lecture Notes in Math., Springer-Verlag, 1996.


\bibitem[{\bf CaOr}]{cascanteortega}
C. Cascante and J.M. Ortega, {\it Norm
inequalities for potential-type operators in homogeneous spaces},
Math. Nachr. {\bf 228} (2001),  85--107.


\bibitem[{\bf CaOrVe1}]{cascanteortegaverbitsky}
C. Cascante, J.M. Ortega and I.E. Verbitsky, {\it Trace
inequalities of Sobolev type in the upper triangle case},
  Proc. London Math. Soc. {\bf 80} (2000), 391--414.

\bibitem[{\bf CaOrVe2}]{cascanteortegaverbitsky2}
C. Cascante, J.M. Ortega and I.E. Verbitsky, 
{\it Wolff's inequality for radially nonincreasing kernels and 
applications to trace inequalities},   Potential Analysis  
{\bf 16} (2002), 347--372;    Erratum,  Potential Analysis   
{\bf 17} (2002), 303--305.


\bibitem[{\bf FeSt}]{feffermanstein}
C. Fefferman and E.M. Stein, {\it Some maximal inequalities},
Amer. J. Math. {\bf 93} (1971), 107--115.
 
\bibitem[{\bf HeWo}]{hedbergwolff}
L.I. Hedberg and Th. H. Wolff, {\it Thin sets in nonlinear
potential theory},
Ann. Inst. Fourier (Grenoble)  
{\bf 33} (1983), 161--187.

 
\bibitem[{\bf KeSa}]{kermansawyer}
R. Kerman  and E. Sawyer, {\it The trace inequality and eigenvalue
estimates for Schr\"odinger operators},
 Ann. Inst. Fourier (Grenoble)  
{\bf 36} (1986), 207--228.


\bibitem[{\bf Ma}]{mazya}
V.G. Maz'ya, {\it Sobolev Spaces}, Springer-Verlag,
Berlin--Heidelberg--New York, 1985.

\bibitem[{\bf MaNe}]{mazyanetrusov}
V.G. Maz'ya  and Y. Netrusov, {\it Some counterexamples for the
theory of Sobolev spaces on bad domains},
 Potential Analysis 
{\bf 4} (1995),  47--65.

 \bibitem[{\bf MaVe}]  {mazyaverbitsky} V.G. Maz'ya and I.E. Verbitsky, 
{\it Capacitary inequalities for
fractional integrals, with applications to partial differential
equations and Sobolev multipliers},
 Arkiv f\"or Matem. {\bf 33} (1995), 81--115.


\bibitem[{\bf NaTrVo}]  {nazarovtreilvolberg} F. Nazarov, S. Treil  
and A. Volberg,  {\it The Bellman functions and two-weight 
inequalities for 
Haar multipliers}, J. Amer. Math. Soc. {\bf 12} (1999), 909--928.


\bibitem[{\bf Sa}]{sawyer} E.T. Sawyer,
{\it A characterization of a two weight norm inequality for
maximal operators}, Studia Math.  {\bf 75} (1982), 1--11.

\bibitem[{\bf SaWh}]{sawyerwheeden}
E.T. Sawyer and R.L. Wheeden, {\it Weighted inequalities for
fractional integrals on Euclidean and homogeneous spaces}, Amer.
J. Math. {\bf 114} (1992), 813--874.



\bibitem[{\bf SaWhZh}] {sawyerwheedenzhao}E.T. Sawyer, R.L.
Wheeden and S. Zhao, {\it Weighted norm inequalities for operators
of potential type and fractional maximal functions}, Potential
Analysis {\bf 5} (1996), 523--580.
 


\bibitem[{\bf Ve1}]{verbitsky2}
I.E. Verbitsky,  {\it Imbedding and multiplier theorems for discrete 
Littlewood--Paley spaces}, 
Pacific J. Math. {\bf 176} (1996), 529--556.

\bibitem[{\bf Ve2}]{verbitsky}
I.E. Verbitsky,  {\it Nonlinear potentials and trace inequalities}, 
Operator Theory: Advances and Appl. {\bf 110} (1999), 323--343.


\bibitem[{\bf VeWh}]{verbitskywheeden}
I.E. Verbitsky and R.L. Wheeden, {\it Weighted norm inequalities
for integral operators}, Trans. Amer. Math. Soc. {\bf 350}
(1998), 3371--3391.







\end{thebibliography}
\end{document}